\documentclass{article}
\usepackage[utf8]{inputenc}
\usepackage[T1]{fontenc}

\usepackage{amsmath,amssymb,amsthm,mathtools}
\usepackage[margin=1.05in]{geometry}
\usepackage{enumitem}
\usepackage[expansion=false]{microtype}
\usepackage{xcolor}
\definecolor{linkc}{rgb}{0,0,0.45}
\usepackage[colorlinks=true,linkcolor=linkc,citecolor=linkc,urlcolor=linkc]{hyperref}
\usepackage[noabbrev]{cleveref}
\AtBeginDocument{\RemoveFromHook{label}[firstaid/cleveref]}

\newtheorem{theorem}{Theorem}[section]
\newtheorem{lemma}[theorem]{Lemma}
\newtheorem{proposition}[theorem]{Proposition}

\newtheorem{assumption}{}

\newtheorem{assumptionB}{}

\newtheorem{assumptionE}{}

\Crefname{assumption}{Assumption}{Assumptions}
\Crefname{assumptionB}{Assumption}{Assumptions}
\Crefname{assumptionE}{Assumption}{Assumptions}
\Crefname{assumptionP}{Assumption}{Assumptions}

\theoremstyle{definition}

\theoremstyle{remark}
\newtheorem{remark}[theorem]{Remark}

\Crefname{conjecture}{Conjecture}{Conjectures}

\newcommand{\R}{\mathbb{R}}
\newcommand{\E}{\mathbb{E}}
\newcommand{\Prob}{\mathbb{P}}
\newcommand{\N}{\mathbb{N}}
\newcommand{\CR}{\mathrm{CR}}
\newcommand{\ind}{\mathbf{1}}

\newcommand{\Fcal}{\mathcal{F}}

\newcommand{\Ecal}{\mathcal{E}}
\newcommand{\Lcal}{\mathcal{L}}
\newcommand{\Zcal}{\mathsf{Z}}
\newcommand{\Ucal}{\mathsf{U}}
\newcommand{\Dcal}{\boldsymbol{\Delta}}
\newcommand{\norm}[1]{\lVert #1\rVert}
\newcommand{\abs}[1]{\lvert #1\rvert}
\newcommand{\TV}{\mathrm{TV}}

\newcommand{\dd}{\mathrm{d}}
\newcommand{\Sym}{\mathfrak{S}}

\title{Linear convergence of evolution strategies to local minima by coupling with quadratic models}

\author{Armand Gissler\thanks{Inria and DI ENS, Paris, France}}

\begin{document}
	\maketitle

	\begin{abstract}
		Among derivative-free optimization methods, evolution strategies have emerged as competitive algorithms for multimodal, noisy or ill-conditioned problems in small and moderate dimensions.
		Theoretical guarantees have long been scarce, but a series of recent works established their convergence on certain classes of objective functions.
		A successful approach normalizes the states of the algorithm into a Markov chain and proves that this chain converges to a stationary distribution.
		It has yielded the linear convergence of several evolution strategies, including the state-of-the-art CMA-ES.
		Yet this approach requires the objective function to be scaling-invariant, hence unimodal, since otherwise the normalized process is not a Markov chain.
		We propose a methodology to prove linear convergence on a much larger class of functions, including smooth multimodal functions.
		Started close enough to a local minimum with a nondegenerate Hessian and with a small enough step size, the $(\mu/\mu_w,\lambda)$-ES with step-size adaptation converges linearly to this minimum with high probability, at the rate of its quadratic model, and the probability of failure vanishes when the initial step size and the initial distance to the minimum tend to zero.
		The proof views the normalized process as an adaptive chain and couples it, with high probability, with the Markov chain associated with the quadratic model.
	\end{abstract}

	\tableofcontents

	\section{Introduction}

	We consider the
	problem of minimization of a numerical function $f\colon\R^d\to\R$.
	In derivative-free optimization, we often assume that 
	$f$ can only be accessed through a black-box oracle without derivatives.
	This setting is addressed by several classes of algorithms: direct search methods \cite{nelder1965simplex,torczon1997convergence,kolda2003optimization,audet2006mesh}; model-based methods \cite{powell2006newuoa,powell2009bobyqa}; randomized methods \cite{nesterov2017random}.

	The focus of this paper is on Evolution Strategies (ES) \cite{rechenberg1973evolutionsstrategie,schwefel1981numerical},
	which are randomized derivative-free algorithms that model a population of candidate solutions in the search space by a probability distribution which ideally tends to a minimum of $f$.
	 The state-of-the-art ES is the Covariance Matrix Adaptation Evolution Strategy (CMA-ES) \cite{hansen2001completely,hansen2003reducing}, which 
	 adapts the parameters of a multivariate normal distribution, including a full covariance matrix, to minimize $f$.
	The ES are particularly competitive in small and moderate dimensions on multimodal, non-separable or noisy problems, for which methods based on local models struggle.

	The theoretical study of ES remains behind their practical use.
	Linear convergence was first established on the sphere function for the elitist $(1+1)$-ES with one-fifth success rule \cite{jagerskupper2007algorithmic}. 
	Global convergence of ES has also been obtained by assuming a sufficent decrease condition~\cite{diouane2015globally}.
	The $(1,\lambda)$-ES with self-adaptation was also proven to be linearly convergent by analyzing the ergodicity of a suitably normalized Markov chain \cite{auger2005convergence}.
	From the latter work,
	in order to prove convergence of ES on scaling-invariant functions~\cite{toure2021scaling},
	 a Markov-chain methodology was developed \cite{auger2016linear}.
	 To simplify such analysis, sufficient conditions for the irreducibility and aperiodicity of such complex Markov chains have been established in recent works~\cite{chotard2019verifiable,gissler2024convergence}.
	 This methodology was then applied to prove linear convergence of step-size adaptive ES on smooth scaling-invariant functions~\cite{toure2023global} and then 
	 of CMA-ES on ellipsoidal functions~\cite{gissler2024linear}.
	 Alternative methods in recent works successfully prove linear convergence of (1+1)-ES on a class of strongly convex functions that may not be scaling-invariant but are still unimodal~\cite{akimoto2022global,morinaga2024linear}.

	The Markov chain approach has an important limitation: the functions on which linear convergence is proven are scaling-invariant hence unimodal~\cite{toure2021scaling}, since the normalized process is a Markov chain only on this class.
	This is far from the situations in which ES are used in practice, where the objective typically has several local minima.
	For such functions, the natural question is a \emph{local} one: if the algorithm is started close enough to a local minimum $x^*$ with a small enough step size, does it converge linearly to $x^*$, and with which probability?
	Note that one cannot hope for an almost-sure statement here: the Gaussian sampling gives a positive probability of leaving any bounded neighbourhood of $x^*$, and once the algorithm has left the area of attraction to $x^*$ nothing can be concluded in general.

	\paragraph{Contributions}
	We answer this question for the $(\mu/\mu_w,\lambda)$-ES with step-size adaptation without cumulation, described in \Cref{sec:algo}.
	Our main result on the algorithm, \Cref{thm:es2}, states that if $f$ is sufficiently regular around a local minimum $x^*$, then for every $\varepsilon\in(0,1)$ there are $r_\varepsilon,s_\varepsilon>0$ such that, whenever $\abs{m_0-x^*}\le r_\varepsilon$ and $\sigma_0\le s_\varepsilon$, with probability at least $1-\varepsilon$ the algorithm converges linearly to $x^*$:
	\[
	\lim_{t\to\infty}\frac1t\log\abs{m_t-x^*}=\lim_{t\to\infty}\frac1t\log\sigma_t=-\CR<0 \enspace,
	\]
	where $\CR$ is the convergence rate of the algorithm on the quadratic model $x\mapsto(x-x^*)^\top\nabla^2f(x^*)(x-x^*)$.
	Nothing is assumed on $f$ outside a neighbourhood of $x^*$.
	In particular $f$ may have any number of other local minima, and $x^*$ need not be the global minimum.
	We obtain a quantitative probability bound (\Cref{thm:es-coupled,th:decoupled-initial-conditions}) which decays as a power of $\sigma_0$.
	Along the way we prove in \Cref{thm:CRquad} that the rate $\CR$ is positive on every ellipsoidal function, under the standard condition that the step size increases on linear functions and for well-chosen hyperparameters.

\begin{figure}[t]
		\centering
		\includegraphics[width=0.8\textwidth]{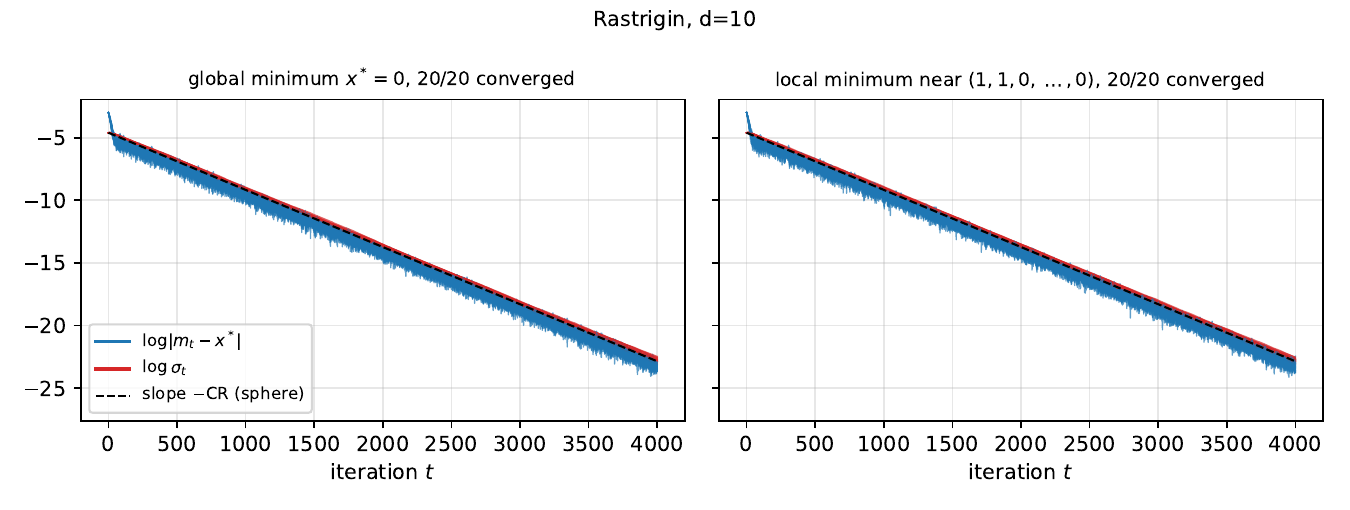}
		\vspace*{-1em}
		\caption{Twenty runs on the Rastrigin function in dimension $d=10$, started at distance $0.05$ of a local minimum $x^*$ with $\sigma_0=0.01$. Left: the global minimum $x^*=0$. Right: the local minimum near $(1,1,0,\dots,0)$. The dashed line has slope $-\CR$ the rate of the quadratic model, estimated by a Monte Carlo run of the normalized chain on the sphere function.}
		\label{fig:exp1}
	\end{figure}

\begin{figure}[t]
		\centering
		\includegraphics[width=0.6\textwidth]{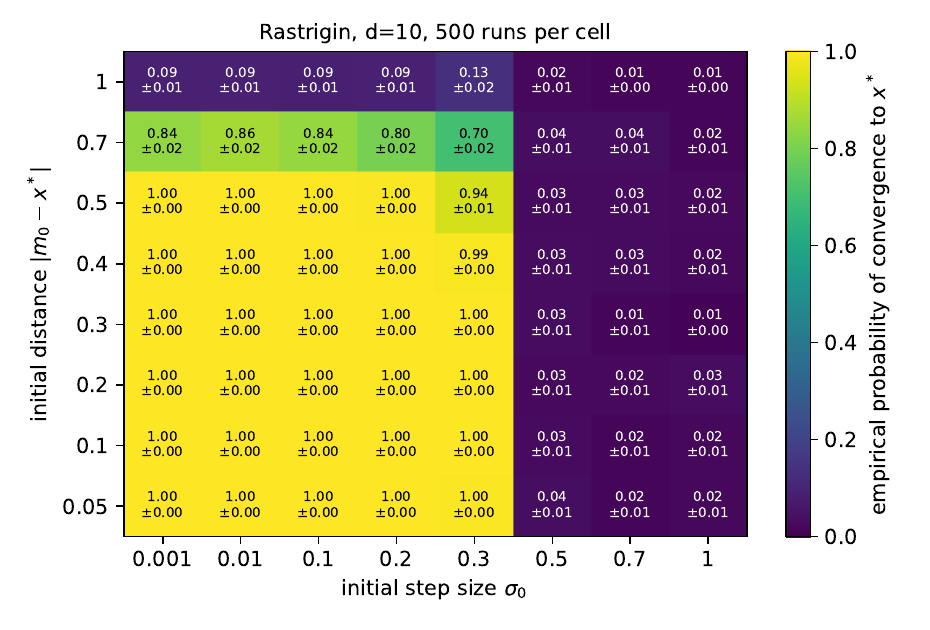}
		\vspace*{-1em}
		\caption{Empirical probability, over $500$ independent runs per cell, that the algorithm started at $m_0=x^*+\abs{m_0-x^*}e$, with $e$ uniform on the unit sphere, and with initial step size $\sigma_0$, converges to the global minimum $x^*=0$ of the Rastrigin function in dimension $d=10$. 
				The $\pm$ values are binomial standard errors.}
		\label{fig:exp2}
	\end{figure}

	\Cref{fig:exp1,fig:exp2} illustrate these results on the Rastrigin function\footnote{We use $f(x)=10d+\sum_{i=1}^d(x_i^2-10\cos(2\pi x_i))$ in dimension $d=10$, with hyperparameters $\lambda=10$, $\mu=5$, the weights $w_i\propto\log(\mu+1)-\log i$, $c_m=0.2$ and $d_\sigma=50$.}.
	Started close to a local minimum $x^*$ with a small step size, the algorithm converges linearly to it, at the rate of its quadratic model.
	The probability of convergence is close to one  as long as the initial step size stays small and $m_0$ is near $x^*$.

	The proof rests on a new point of view on the normalized process, which we believe to be of independent interest.
	On a general objective function, the normalized state $z_t=(m_t-x^*)/\sigma_t$ is no longer a Markov chain, because the ranking of the candidate solutions by $f$ depends on $\sigma_t$ through the function $x\mapsto f(x^*+\sigma_tx)$.
	It is, however, an \emph{adaptive chain}: its transition kernel $P_\sigma$ belongs to a family indexed by the step size $\sigma$, and the index $\sigma_t$ is itself updated along the trajectory.
	As $\sigma\to0$, the function $x\mapsto f(x^*+\sigma x)$ seen by the normalized process approaches, up to an increasing transformation, its quadratic model, which is scaling-invariant.
	The kernels $P_\sigma$ therefore approach the limit kernel $P_0$ of the normalized chain on the quadratic model, whose geometric ergodicity is known from \cite{toure2023global}.
	In \Cref{sec:adaptive} we develop this idea in an abstract setting.
	\Cref{thm:main} couples the adaptive chain with the limit chain driven by $P_0$ and the same innovations: outside an event of probability at most $CV(z_0)\sigma_0^\theta$, where $V$ is the Lyapunov function of $P_0$, the two chains coincide for all times, so that the adaptive chain inherits the almost-sure limit theorems of the limit chain.
	The proof combines a perturbation argument for the exponential moments of the step size, based on the spectral theory of quasi-compact operators \cite{hennion2001limit,kontoyiannis2003spectral}, with a quantitative bound on the probability that the two kernels disagree, which for the ES follows from a Taylor expansion and an anti-concentration inequality for Gaussian polynomials \cite{carbery2001distributional}.
	\Cref{thm:main} alone requires $\abs{m_0-x^*}$ to be small compared to $\sigma_0$.
	\Cref{th:decoupled-initial-conditions} removes this restriction by coupling only after a random time at which the normalized state has entered a bounded set, the existence and the cost of which are controlled by the drift analysis of H\'ajek \cite{hajek1982hitting}, a tool that is standard in the theory of evolutionary algorithms on discrete search spaces \cite{he2001drift,lengler2020drift}.

	\paragraph{Related works on adaptive chains}
	Processes whose transition kernel depends on a parameter updated along the trajectory are the controlled Markov chains of stochastic approximation theory \cite{benveniste1990adaptive,kushner2003stochastic}
	or the Markov decision processes in reinforcement learning \cite{konda2003actor}.
	Under the name of adaptive Markov chains, they have been studied mostly in the context of adaptive MCMC, see \cite{haario2001adaptive,andrieu2006ergodicity,roberts2007coupling,andrieu2008tutorial}.
	There, the transition kernels $P_\theta$ usually share a common target distribution $\pi$, and the parameter $\theta_n$ is tuned by the user to speed up the convergence to $\pi$.
	The classical conditions are diminishing adaptation, which asks that $P_{\theta_{n+1}}$ and $P_{\theta_n}$ become close, and bounded convergence, which asks that the frozen kernels $P_{\theta_n}$ converge to $\pi$ in a time that stays bounded in probability \cite{roberts2007coupling}.
	The latter is also called containment~\cite{fort2011convergence}, and it follows from simultaneous uniform ergodicity \cite{roberts2007coupling} or from a drift and a minorization condition with constants uniform in $\theta$ \cite{andrieu2006ergodicity}.
	 This framework has also been studied with targets $\pi_{\theta}$ which depend on $\theta$,  with convergence of the marginals and a strong law of large numbers under uniform geometric ergodicity of the kernels, diminishing adaptation in $V$-norm and convergence of the targets \cite{fort2011convergence}.
	Recent work~\cite{chotard2026ergodicity} proposed an adaptive Metropolis--Hastings sampler whose step-size and covariance adaptation are borrowed from CMA-ES.

	Our setting differs from these frameworks in several respects.
	The kernels $P_\sigma$ for $\sigma>0$ have no invariant distribution in general and
	only the limit kernel $P_0$ is assumed ergodic.
	The parameter converges to the boundary point $\sigma=0$ of the parameter space, where the kernels degenerate.
	Diminishing adaptation would mean here that the step size converges, which is the conclusion we are after, so it is not assumed but proven.
	In its place we assume a quantitative convergence of $P_\sigma$ to $P_0$ as $\sigma\to0$, together with a sign condition on the expectation of $\log\sigma_t$ under the invariant distribution of $P_0$.
	The multiplicative form of the adaptation requires exponential moments where polynomial moments suffice in \cite{fort2011convergence}.
	The role of containment is played by a drift condition in the form of Foster and H\'ajek towards the region where the coupling is available.
	Finally, instead of the convergence of the marginals, we obtain that the adaptive chain coincides with the limit chain on an event of high probability, which yields almost-sure limits and the rate of decay of the parameter.
	To our knowledge, this is the first use of adaptive chains for the analysis of an ES algorithm.

	\paragraph{Scope and limitations}
	The algorithm analyzed here is a step-size adaptive ES, not CMA-ES: there is no covariance matrix adaptation, and the step-size adaptation uses no cumulation path.
	This setting is usual in the theoretical literature on the Markov-chain methodology \cite{auger2016linear,toure2023global}, and the extension to covariance matrix adaptation, in the spirit of \cite{gissler2024linear}, is left for future work.
	\Cref{thm:CRquad}, and hence \Cref{thm:es2}, is stated for a small learning rate $c_m$ and a large damping $d_\sigma$.
	We believe the rate $\CR$ to be positive for the default parameters as well, as numerical evidence suggests.

	\paragraph{Organization of the paper}
	\Cref{sec:algo} describes the algorithm, recalls the Markov-chain methodology on scaling-invariant functions, and states our results on the algorithm, \Cref{thm:CRquad,thm:es2}.
	\Cref{sec:adaptive} introduces the abstract adaptive chain, its assumptions and the two coupling results, \Cref{thm:main,th:decoupled-initial-conditions}, which are proven in \Cref{sec:proofmain,sec:proofdecoupled}.
	\Cref{sec:proof-es} verifies the assumptions for the ES and proves \Cref{thm:es2}, and \Cref{sec:CRquad} proves \Cref{thm:CRquad}.

	\paragraph{Notation}
	The Euclidean norm on $\R^d$ is denoted by $\abs\cdot$, and we write $\N=\{0,1,\dots\}$, $\R_{++}=(0,\infty)$, and $\abs{u}_\infty:=\max_i\abs{u^i}$ for $u=(u^1,\dots,u^\lambda)\in(\R^d)^\lambda$.
	For a measure $\mu$ and a function $\varphi$ we write $\mu(\varphi):=\int\varphi\dd\mu$, and for a Markov kernel $P$, $P\varphi(z):=\int\varphi(y)P(z,\dd y)$.
	For a measurable $V\ge1$, the $V$-norm of a function is $\norm\varphi_V:=\sup\abs\varphi/V$ and the $V$-norm of a signed measure is $\norm\mu_V:=\sup_{\norm\varphi_V\le1}\abs{\mu(\varphi)}$.
	The total variation norm $\norm{\cdot}_\TV$ is the $V$-norm for $V\equiv1$.
	We use the terminology of \cite{meyn2009markov} for Markov chains on a general state space: $\varphi$-irreducibility, aperiodicity, petite sets, geometric ergodicity.
	We write $\Fcal_t$ for the $\sigma$-field generated by the samples up to iteration $t$, and $\Sym_\lambda$ for the set of permutations of $\{1,\dots,\lambda\}$.

	\section{Local convergence of the $(\mu/\mu_w,\lambda)$-ES}\label{sec:algo}

	\subsection{Description of the algorithm}
	
	Given an objective function $f\colon\R^d\to\R$, the goal of an optimization algorithm is to obtain state variables which form a sequence converging to a minimum $x^*$ of $f$.
	The state variables of  the $(\mu/\mu_w,\lambda)$-ES are $\{(m_t,\sigma_t)\}_{t\in\N}$ where $m_t\in\R^d$  is the mean variable and $\sigma_t>0$ is the step size.
	
	Fix $\lambda\ge2$ an integer called the population size and $\mu\in\{1,\dots,\lambda\}$ the parent number.
	Fix as well positive weights $w=(w_1,\dots,w_\mu)\in\R^\mu_{++}$ that sum to one, a learning rate $c_m>0$ and a damping parameter $d_\sigma>0$.
	Also consider the step-size change $\Gamma\colon\R^d\to\R_{++}$ defined by
	\begin{equation}\label{eq:CSA2}
	\Gamma(y)= \exp\left( \frac{1}{2d_\sigma} \times \left(   \frac{\abs y^2}{d \abs w^2} - 1 \right)  \right) \enspace.
	\end{equation}
	At iteration $t\in\N$, the update is described as follows.
	Consider $U_{t+1}^1,\dots,U_{t+1}^\lambda$ $\lambda$ i.i.d.\ samples from $\mathcal N(0,I_d)$ independent from the past iterations and sort them, ties being broken by the index of the sample, so that
	\begin{equation}
		f( m_t + \sigma_t U_{t+1}^{1:\lambda}) \le \cdots \le f( m_t + \sigma_t U_{t+1}^{\lambda:\lambda}) \enspace.
	\end{equation}
	Define then the mean displacement
	\begin{equation}
		Y_{t+1} := \sum_{i=1}^\mu w_i U^{i:\lambda}_{t+1}
	\end{equation}
	and set
	\begin{equation}\label{eq:es}
		m_{t+1}:=m_t+ c_m \sigma_tY_{t+1},\qquad \sigma_{t+1}:=\sigma_t\Gamma(Y_{t+1}) \enspace.
	\end{equation}

	\subsection{Methodology to prove linear convergence on scaling-invariant functions}

	To prove linear convergence of ES, \cite{auger2016linear,toure2023global} propose to study the normalized process $\{z_t\}_{t\in\N}$ defined by
	\begin{equation}
		z_t = \frac{m_t-x^*}{\sigma_t} \enspace,
	\end{equation}
	and to show that it converges to a stationary distribution $\pi$.
	The step-size recursion in \eqref{eq:es} gives
	\begin{equation}\label{eq:LLN}
		\frac{1}{t} \log \abs{m_t-x^*} = \frac{1}{t} \log\sigma_0 + \frac{1}{t} \sum_{k=0}^{t-1} \log \Gamma (Y_{k+1}) + \frac1t\log\abs{z_t} \enspace.
	\end{equation}
	If $\frac1t\log\abs{z_t}\to0$ and a law of large numbers applies to the sum, the right-hand side converges to
	\begin{equation}\label{eq:CRdef}
		-\CR:=\E_\pi\big[\log\Gamma(Y_1)\big] \enspace,
	\end{equation}
	where $\E_\pi$ refers to the process started from $z_0\sim\pi$.
	Suppose that $f$ is scaling-invariant with respect to $x^*$, that is~(see~\cite{toure2021scaling}):
	\begin{equation}
		f(x^* + x) \le f(x^*+y) \iff  f(x^*+a x) \le f(x^*+a y) , \qquad \forall x,y\in\R^d, \forall a>0 \enspace.
	\end{equation}
	The ranking of $m_t+\sigma_tU^1_{t+1},\dots,m_t+\sigma_tU^\lambda_{t+1}$ by $f$ is then the one of $z_t+U^1_{t+1},\dots,z_t+U^\lambda_{t+1}$ by $f(x^*+\cdot)$.
	Hence $Y_{t+1}$ is a function of $z_t$ and $\mathbf U_{t+1}$, and \eqref{eq:es} gives
	\begin{equation}\label{eq:zrec}
		z_{t+1}=\frac{z_t+c_mY_{t+1}}{\Gamma(Y_{t+1})} \enspace.
	\end{equation}
	The process $\{z_t\}_{t\in\N}$ is therefore a time-homogeneous Markov chain \cite[Proposition 2]{toure2023global}.

	Let $\mathcal N_{1:\lambda}\le\dots\le\mathcal N_{\lambda:\lambda}$ be the order statistics of $\lambda$ i.i.d.\ standard normal variables.
	On a linear function $x\mapsto e^\top x$ with $\abs e=1$, the ranking of the samples is the ranking of their coordinates along $e$, and the coordinates orthogonal to $e$ are not affected by the selection.
	By \eqref{eq:CSA2}, this gives
	\begin{equation}\label{eq:CRlin}
		\E\big[\log\Gamma(Y_{t+1})\big]=\frac{1}{2d_\sigma d\abs w^2}\Big(\E\Big[\Big(\sum_{i=1}^\mu w_i\mathcal N_{i:\lambda}\Big)^2\Big]-\abs w^2\Big) \enspace,
	\end{equation}
	which depends neither on $e$ nor on $(m_t,\sigma_t)$.
	The step size increases on average on linear functions when this quantity is positive, that is, when
	\begin{equation}\label{eq:lin-div}
		\E\Big[\Big(\sum_{i=1}^\mu w_i\mathcal N_{i:\lambda}\Big)^2\Big]>\abs w^2 \enspace.
	\end{equation}

	\begin{proposition}[{\cite[Lemma 2, Propositions 7, 11 and 12 and Theorems 6 and 7]{toure2023global}}]\label{prop:TAH}
		Let $f$ be continuously differentiable and scaling-invariant with respect to $x^*$, and let $x^*$ be its unique global minimizer.
		Assume \eqref{eq:lin-div}.
		Then $\{z_t\}_{t\in\N}$ is a $\varphi$-irreducible and aperiodic Markov chain, every compact set is petite, and there is $\beta_0>0$ such that, for every $\beta\in(0,\beta_0]$, there are $\lambda_{\beta}\in(0,1)$, $b_{\beta}<\infty$ and $r_{\beta}>0$ with
		\begin{equation}\label{eq:TAH-drift}
			\E\big[1+\abs{z_{t+1}}^{\beta}\big\vert z_t=z\big]\le\lambda_{\beta}\big(1+\abs z^{\beta}\big)+b_{\beta}\ind_{\{\abs z\le r_{\beta}\}} \qquad\text{for all }z\in\R^d \enspace.
		\end{equation}
		In particular $\{z_t\}_{t\in\N}$ is geometrically ergodic with a unique invariant probability $\pi$, and $\log\abs\cdot$ is $\pi$-integrable.
		Moreover, for every initial condition with $m_0\ne x^*$, almost surely,
		\begin{equation}\label{eq:TAH-LLN}
			\lim_{t\to\infty}\frac1t\log\abs{m_t-x^*}=\lim_{t\to\infty}\frac1t\log\sigma_t=-\CR \enspace.
		\end{equation}
	\end{proposition}

	By \eqref{eq:TAH-LLN}, the algorithm converges linearly to $x^*$ if and only if $\CR>0$.
	The sign of $\CR$ is not given by \Cref{prop:TAH}.
	The following theorem gives it on quadratic functions, for a large damping factor and a small learning rate.

	\begin{theorem}[Positive convergence rate on quadratic functions]\label{thm:CRquad}
		Let $f(x)=(x-x^*)^\top H(x-x^*)$ with $H$ symmetric positive definite.
		Assume \eqref{eq:lin-div}, that the weights are nonincreasing, $w_1\ge\dots\ge w_\mu$, and that $\mu\le\lambda/2$.
		There are $\bar c>0$ and $\bar d<\infty$ such that $\CR>0$ when $c_m\le\bar c$ and $c_md_\sigma\ge\bar d$.
	\end{theorem}
	The proof is given in \Cref{sec:CRquad}.

	\subsection{Assumptions and main result}
	
	Our first assumption is on the objective function $f$. Note that ES are invariant to increasing transformations, so any objective function $h\circ f$ with $h\colon\R\to\R$ increasing and $f$ satisfying \ref{ass:function}
	would still satisfy our convergence result below.

	\begin{assumptionE}\label{ass:function}
		We suppose that $f$ is measurable on $\R^d$ and twice continuously differentiable on a ball $B(x^*,r_0)$.
		Besides, we suppose $\nabla f(x^*)=0$, $\nabla^2f(x^*)\succ 0$ and $\nabla^2 f$ is locally Lipschitz on $B(x^*,r_0)$. 
	\end{assumptionE}
	
	Second, we consider the normalized process $\{z_t\}_{t\in\N}$ defined by
	\begin{equation}\label{eq:normalized-chain}
		z_t = \frac{m_t-x^*}{\sigma_t} \enspace.
	\end{equation}
	Let $f_0(x):=(x-x^*)^\top\nabla^2f(x^*)(x-x^*)$ be the quadratic model of $f$ at $x^*$.
	It is scaling-invariant with respect to $x^*$, so \Cref{prop:TAH} and \Cref{thm:CRquad} apply to it under the following assumption on the parameters of the algorithm.

	\begin{assumptionE}\label{ass:quad}
		The weights are nonincreasing, $w_1\ge\dots\ge w_\mu$, the parent number satisfies $\mu\le\lambda/2$, and \eqref{eq:lin-div} holds.
		Moreover $c_m\le\bar c$ and $c_md_\sigma\ge\bar d$, where $\bar c$ and $\bar d$ are the constants of \Cref{thm:CRquad} for $H=\nabla^2f(x^*)$.
	\end{assumptionE}

	Under \ref{ass:function} and \ref{ass:quad}, \Cref{prop:TAH} and \Cref{thm:CRquad} apply to $f_0$.
	The normalized process of the algorithm \eqref{eq:es} on $f_0$ is then a geometrically ergodic Markov chain, which satisfies the drift \eqref{eq:TAH-drift} with the potential
	\begin{equation}\label{eq:Vbeta}
		V_\beta(z):=1+\abs z^\beta
	\end{equation}
	for every $\beta\in(0,\beta_0]$, and its convergence rate $\CR$ defined via \eqref{eq:CRdef} is positive.
	We fix $\beta\in(0,\min(1,\beta_0)]$.
	The following convergence result is obtained by comparing the normalized process of the algorithm on $f$ with the one on $f_0$.
	\Cref{sec:proof-es} is devoted to its proof.

	\begin{theorem}[Local linear convergence of the $(\mu/\mu_w,\lambda)$-ES]\label{thm:es2}
		Assume \ref{ass:function} and \ref{ass:quad}. For every $\varepsilon\in(0,1)$ there exist $r_\varepsilon>0$ and $s_\varepsilon>0$ such that, whenever
		\[
		\abs{m_0-x^*}\le r_\varepsilon\qquad\text{and}\qquad0<\sigma_0\le s_\varepsilon\enspace,
		\]
		there is an event of probability $\ge1-\varepsilon$ on which
		\[
		\lim_{t\to\infty}\frac1t\log\sigma_t=\lim_{t\to\infty}\frac1t\log\abs{m_t-x^*}=-\CR \enspace,
		\]
		where $\CR>0$ is the same convergence rate as for quadratic objective functions with Hessian matrices proportional to $\nabla^2 f(x^*)$.
	\end{theorem}

	\section{Convergence of an adaptive chain}\label{sec:adaptive}

	\subsection{Coupling with a limit Markov chain}
	
	Consider the stochastic process $\{z_{t}\}_{t\in\N}$ 
	obeying the recursion
	\begin{equation}\label{eq:adaptive-z}
		z_{t+1} = \Phi_{\sigma_{t}}(z_t,\mathbf U_{t+1})
	\end{equation}
	where $\Phi \colon  (\sigma,z,u) \in (0,+\infty)\times\Zcal\times \Ucal   \mapsto
	\Phi_\sigma (z,u) \in\Zcal$ is a measurable function, $\sigma_t>0$ is an adaptive parameter, and
	$\mathbf U = \{\mathbf U_{t+1}\}_{t\in\N}$ is an i.i.d.\ process with common probability distribution $\nu$ on a measurable space $\Ucal$.
	Defining $P_\sigma(\cdot|z)$ as the probability distribution of $\Phi_\sigma (z,U_1)$, we can write the transition law of $\{z_t\}_{t\in\N}$ as
	\begin{equation}
		z_{t+1} \sim P_{\sigma_{t}} (\cdot | z_t ) \enspace.
	\end{equation}
	With the observation that $\{P_\sigma\}_{\sigma>0}$ is a family of Markov transition kernels on $\Zcal$ and  $\sigma_t>0$ is an adaptive parameter, we interpret $\{z_t\}_{t\in\N}$ as an adaptive chain.

	The adaptive chain is compared with a limit chain.
	Consider the \emph{limit map} $\Phi_0\colon\Zcal\times\Ucal\to\Zcal$ to be measurable, and let $P_0(\cdot|z)$ be the probability distribution of $\Phi_0(z,\mathbf U_1)$.
	Note that $\Phi_0$ is not simply the point-wise limit of $\Phi_\sigma$ and the value $\sigma=0$ is not reached by the adaptive parameter.
	The two are related by \ref{ass:coupling} below.
	In all statements, the case $\sigma=0$ of a property required for every $\sigma\ge0$ refers to the limit map.
	For the ES, $\Phi_\sigma$ is the map of the algorithm on $f$ with step size $\sigma$, and $\Phi_0$ is the map of the algorithm on the quadratic model $f_0$, see \Cref{sec:es-notation}.

	The following assumption defines the adaptation rules of $\sigma_t$ that we consider.

	\begin{assumption}\label{ass:normalized}
		There exists a measurable map $\Gamma: (0,+\infty) \times\mathsf Z\times\mathsf U\to(0,\infty)$, such that for every $t\ge0$,
		\begin{equation}\label{eq:normalized-dynamics-sigma}
			\sigma_{t+1}=\sigma_t\times \Gamma_{\sigma_t}(z_t,\mathbf U_{t+1}) \enspace.
		\end{equation}
		There is also a measurable limit rule $\Gamma_0\colon\Zcal\times\Ucal\to(0,\infty)$, attached to the limit map $\Phi_0$.
		Besides, there is a function $\ell\colon\Ucal\to(0,+\infty)$ such that $\abs{\log \Gamma_\sigma (z,u)} \le \ell (u)$ for every $u\in\Ucal$, $\sigma\ge0$ and $z\in\Zcal$.
	\end{assumption}

	The limit chain $\{(z^0_t,\sigma^0_t)\}_{t\in\N}$ is the process driven by the limit map and the limit rule with the same i.i.d.\ random process $\mathbf U$,
	\begin{equation}\label{eq:limit-chain}
		z^0_{t+1}=\Phi_0(z^0_t,\mathbf U_{t+1}),\qquad \sigma^0_{t+1}=\sigma^0_t\times\Gamma_0(z^0_t,\mathbf U_{t+1}) \enspace,
	\end{equation}
	started from $(z^0_0,\sigma^0_0)=(z_0,\sigma_0)$.
	The process $\{z^0_t\}_{t\in\N}$ is a Markov chain with transition kernel $P_0$, and $\log\sigma^0_t$ is an additive functional of its path.
	We also use the limit chain started at a time $s\in\N$ from a state $(z,\sigma)$, which is defined by \eqref{eq:limit-chain} for $t\ge s$ with $(z^0_s,\sigma^0_s)=(z,\sigma)$.
	For the ES, the limit chain is the normalized process and the step size of the algorithm on the quadratic model $f_0$.

	As a replacement of the uniform ergodicity assumption, we only assume the ergodicity of the limit kernel $P_0$, via the following drift condition.
	
	\begin{assumption}\label{ass:potential}\label{ass:erg}
		The transition kernel $P_0$ is $\varphi$-irreducible and aperiodic, and there are
		 a measurable potential function $V:\Zcal\to[1,\infty)$,
		 $\lambda_V\in(0,1)$, $b_V<\infty$ and a petite set $\mathsf C$ such that
		\[
		P_0V\le\lambda_V V+b_V\ind_{\mathsf C} \enspace.
		\]
		Besides, there exists a measurable function $\Xi:\Ucal\to[1,\infty)$
		such that
		\[
		V\big(\Phi_\sigma(z,u)\big)\le V(z)\Xi(u)
		\]
		for all $u\in\Ucal, \sigma\ge0, z\in\Zcal$.
	\end{assumption}

	By \cite[Theorem 15.0.1]{meyn2009markov}, \Cref{ass:erg} implies that $P_0$ has a unique invariant probability $\pi_0$, that $\pi_0(V)<\infty$, and that there are $M<\infty$ and $\rho_V\in(0,1)$ with
	\begin{equation}\label{eq:Vgeo}
		\norm{P_0^n(z,\cdot)-\pi_0}_V\le M\rho_V^nV(z)   \enspace.
	\end{equation}
	Yet, most convergence results on adaptive chains assume diminishing adaptation~\cite{andrieu2006ergodicity,roberts2007coupling,fort2011convergence}, which would translate here into a convergence of $\sigma_t$ to $0$. 
	We replace it here by the following two assumptions, which ensure the coupling of the adaptive chain $\{z_t\}_{t\in\N}$ with the limit chain \eqref{eq:limit-chain}.

	\begin{assumption}\label{ass:coupling}
		There exist $\bar\sigma>0$, $\kappa<\infty$,  $\alpha>0$ and $\eta>0$ such that, with
		\[
		\Dcal_\sigma(z):=\left\{u\in\Ucal:\ \Phi_\sigma(z,u)\ne\Phi_0(z,u)\ \text{or}\ \Gamma_\sigma(z,u)\ne \Gamma_0(z,u)\right\},
		\]
		one has $\nu(\Dcal_\sigma(z))\le\kappa\sigma^\eta$ for all $\sigma\le\bar\sigma$ and all $z\in\Zcal_\sigma$, where
		\begin{equation}\label{eq:Zsigma}
			\Zcal_\sigma:=\big\{z\in\Zcal \mid V(z)\le\sigma^{-\alpha}\big\} \qquad \text{ for } \sigma>0 \enspace.
		\end{equation}
	\end{assumption}
	
	We set for $\sigma\ge0$:
	\begin{equation}\label{eq:def-P0}
		\gamma_\sigma(z):=\int\log \Gamma_\sigma(z,u)\nu(\dd u) \enspace.
	\end{equation}
	
	\begin{assumption}\label{ass:rate}
		We have	$\CR:=-\pi_0(\gamma_0)>0$.
	\end{assumption}
	
	Finally, we require the following integration assumption.
	
	\begin{assumption}\label{ass:domination}
		There exists $\Lambda>0$ such that
		the functions $\Xi(\cdot)$ and $\ell(\cdot)$ from \ref{ass:erg} and \ref{ass:normalized} satisfy
		$\int e^{\Lambda\ell(\cdot)}\Xi (\cdot) \dd\nu<\infty$.
	\end{assumption}

	We have now introduced all the mathematical background to state the following coupling result.
	
	\begin{theorem}\label{thm:main}
		Assume \ref{ass:normalized}--\ref{ass:domination}.
		There are $\theta>0$ and $C<\infty$ such that, for every initial condition $(z_0,\sigma_0)\in\Zcal\times\R_{++}$, there is an event $\Ecal$ with
		\begin{equation}\label{eq:B}
			\Prob(\Ecal^c)\le CV(z_0)\sigma_0^{\theta}
		\end{equation}
		on which, for every $t\ge0$, $(z_t,\sigma_t)=(z^0_t,\sigma^0_t)$, where $\{(z^0_t,\sigma^0_t)\}_{t\in\N}$ is the limit chain \eqref{eq:limit-chain} started from $(z_0,\sigma_0)$. Consequently, up to a null set, on $\Ecal$
		\[
		\lim_{t\to\infty}\frac1t\log\sigma_t=-\CR
		\qquad\text{and}\qquad
		\lim_{t\to\infty}\frac1t\log V(z_t)=0 \enspace,
		\]
		and, for all $t\ge0$,
		\[
		\norm{\Lcal(z_t)-\pi_0}_\TV\le
		2CV(z_0)\sigma_0^{\theta}
		+\norm{P_0^t(z_0,\cdot)-\pi_0}_\TV \enspace.
		\]
	\end{theorem}
	\Cref{thm:main} is proven through \Cref{sec:proofmain}.
	While the bound \eqref{eq:B} allows us to state a convergence result for ES with high probability, it would require initial conditions $\abs{z_0}$ and $\sigma_0$ to be small, resulting in needing the distance of the initial mean to the optimum to be sufficiently smaller than $\sigma_0$, which in practice can be hard to match. 
	The next section is devoted to extend the convergence results to more initial conditions.
	
	\subsection{Coupling after a random time}
	
	In this section, instead of coupling the whole process $\{(z_t,\sigma_t)\}_{t\in\N}$ to a reference Markov process, we introduce a random time from which the coupling occurs.
	This is based on drift theory results~\cite{hajek1982hitting,lengler2020drift}, resulting in the following assumption.
	
	\begin{assumptionB}\label{assB:drift}
		There exist measurable functions $D\colon\Zcal\to\R$
		and $g\colon\Zcal\times\R_{++}\to\R_+$ such that the stopping time
		\begin{equation}\label{eq:B1-def-tau}
			\tau = \inf\{ t\ge 0 \mid D(z_t)\le0 \text{ or } g(z_t,\sigma_t)>1 \}
		\end{equation}
		and some constants $c>0$, $\vartheta>0$, $\varsigma>0$ and $K\ge e^{\vartheta c}$ satisfy for $t < \tau$:
		\begin{equation}\label{eq:drift-for-zsig}
		\begin{array}{ll}
		&	\E[D(z_t)-D(z_{t+1})|\Fcal_t] \ge c \\
		&	\E[e^{\vartheta\abs{ D(z_t) - D(z_{t+1}) }  } | \Fcal_t ] \le K \\
		&	\E[g(z_{t+1},\sigma_{t+1} )  | \Fcal_t  ]   \le g(z_t,\sigma_t)  \\
		&	\sigma_t \le \varsigma\times e^{-D(z_t)}
		\end{array}
		\end{equation}
		Moreover, the potential $V$ of \ref{ass:erg} is bounded on $\{D\le0\}$.
	\end{assumptionB}
	
	The following theorem is our last major contribution. It is the final key to obtain convergence of ES to a minimum $x^*$ of $f$ for well-chosen (but not unreasonable) initial conditions.
	The proof of \Cref{th:decoupled-initial-conditions} can be found in \Cref{sec:proofdecoupled}.
	
	\begin{theorem}\label{th:decoupled-initial-conditions}
		Assume \ref{ass:normalized}--\ref{ass:domination} and \ref{assB:drift}.
		Then, the stopping time $\tau$ from \eqref{eq:B1-def-tau} is almost surely finite. 
		Besides, there are $\theta>0$ and $C<\infty$ such that, for every initial condition $(z_0,\sigma_0)\in\Zcal\times\R_{++}$, there is an event $\Ecal_\tau$ with
		\begin{equation}
			\Prob(\Ecal_\tau^c)  \le g(z_0,\sigma_0)  + C \sup_{(D\le 0)} V \times \left(  \sigma_0^\theta + C_{\rm occ}\varsigma^\theta \int e^{\theta\ell} \dd\nu   \right)
		\end{equation}
		where $C_{\rm occ}<\infty$ depends only on $(c,\vartheta, \theta,K)$, on which, for every $t\ge\tau$, $(z_t,\sigma_t)=(z^0_t,\sigma^0_t)$, where $\{(z^0_t,\sigma^0_t)\}_{t\ge\tau}$ is the limit chain \eqref{eq:limit-chain} started at time $\tau$ from $(z_\tau,\sigma_\tau)$. Consequently, up to a null set, on $\Ecal_\tau$
		\[
		\lim_{t\to\infty}\frac1t\log\sigma_t=-\CR
		\qquad\text{and}\qquad
		\lim_{t\to\infty}\frac1t\log V(z_t)=0 \enspace.
		\]
	\end{theorem}

	\section{Proof of \Cref{thm:main}}\label{sec:proofmain}

		\subsection{Coupling with the limit chain}\label{sec:coupling}
	
	Let
	$\{(z_t,\sigma_t)\}_{t\in\N}$ be the adaptive chain \eqref{eq:adaptive-z} and
	 $\{(z_t^0,\sigma_t^0)\}_{t\in\N}$ be the limit chain \eqref{eq:limit-chain} started from $(z_0,\sigma_0)$.
	
	We define the event $\Ecal$ by 
	\begin{equation}\label{eq:Eevent}
		\Ecal =   \{ \forall t \in\N, ~ \mathbf U_{t+1}\not\in\boldsymbol{\Delta}_{\sigma_t}(z_t) \}
		\cap
		\{ \forall t \in\N , ~ z_t \in \Zcal_{\sigma_t}\}
		\cap
		\{  \forall t \in\N , ~ \sigma_t\le \bar\sigma  \}
		 \enspace,
	\end{equation}
	where the set $\boldsymbol{\Delta}_\sigma(z)$ is defined in \ref{ass:coupling}.
	
	\begin{proposition}
		Assume \ref{ass:normalized}.
		If $\Ecal$ occurs, then for all $t\in\N$, $(z_t^0,\sigma_t^0)=(z_t,\sigma_t)$.
	\end{proposition}
	\begin{proof}
		Immediate with the definitions.
	\end{proof}
	
	\begin{proposition}\label{prop:bound-Ec}
		Assume \ref{ass:normalized} and \ref{ass:coupling}.
		We have 
		\begin{equation}\label{eq:bound-E}
			\Prob(\Ecal^c)\le 
			\kappa\sum_{t=0}^{+\infty}\E[
			\ind_{\{\sigma^0_t\le\bar\sigma\}}(\sigma_t^0)^\eta]
			+ \sum_{t=0}^{+\infty}\Prob [V(z_t^0)>(\sigma_t^0)^{-\alpha}]
			+   \sum_{t=0}^{+\infty} \Prob[ \sigma_t^0>\bar\sigma ] 
		\end{equation}
	\end{proposition}
	
	\begin{proof}
		Let $T$ be the following random time:
		\[
		T:=\inf\big\{t\ge0:\ \mathbf U_{t+1}\in\Dcal_{\sigma_t}(z_t) \text{ or } z_t \not\in \Zcal_{\sigma_t} 
		\text{ or } \sigma_t>\bar \sigma
		\big\} \enspace,
		\]
		so that $\Ecal^c=\{T<\infty\}$ and $\Prob(\Ecal^c)=\sum_{t\ge0}\Prob(T=t)$.
		 By induction, we obtain $(z_t,\sigma_t)=(z^0_t,\sigma^0_t)$ for $t\le T$.
		By conditioning on $\Fcal_t$, we have
		\begin{align*}
		\Prob(T=t\mid \Fcal_t)
		&\le \E\Big[\ind_{\{T\ge t\}}
		\ind_{\{\sigma_t\le\bar\sigma\}}
		\ind_{\{z_t\in\Zcal_{\sigma_t}\}}
		\nu\big(\Dcal_{\sigma_t}(z_t)\big) \mid \Fcal_t \Big] \\
		&\quad+ \Prob [T\ge t \text{ and } z_t\not\in\Zcal_{\sigma_t} \mid \Fcal_t  ]
		+ \Prob[T\ge t \text{ and } \sigma_t>\bar\sigma \mid \Fcal_t ] \\
		&=\E\Big[\ind_{\{T\ge t\}}
		\ind_{\{\sigma_t^0\le\bar\sigma\}}
		\ind_{\{z^0_t\in\Zcal_{\sigma^0_t}\}}\nu\big(\Dcal_{\sigma^0_t}(z^0_t)\big) \mid \Fcal_t  \Big] \\
		&\quad+ \Prob [T\ge t \text{ and } z_t^0\not\in\Zcal_{\sigma_t^0} \mid \Fcal_t  ]
		+ \Prob[T\ge t \text{ and } \sigma_t^0>\bar\sigma \mid \Fcal_t ] \\
		&\le \E\big[
		\ind_{\{\sigma^0_t\le\bar\sigma\}}
		\ind_{\{z_t^0\in\Zcal_{\sigma_t^0}\}}\nu(\Dcal_{\sigma^0_t}(z^0_t)) \mid \Fcal_t  \big] \\
		&\quad+ \Prob [z_t^0\not\in\Zcal_{\sigma_t^0} \mid \Fcal_t  ]
		+ \Prob[ \sigma_t^0>\bar\sigma \mid \Fcal_t ] \\
		&\le
		\kappa(\sigma_t^0)^\eta
		\ind_{\{\sigma_t^0\le\bar\sigma\}}
		+ \Prob [ V(z_t^0)>(\sigma_t^0)^{-\alpha} \mid \Fcal_t  ]
		+ \Prob[ \sigma_t^0>\bar\sigma \mid \Fcal_t ]
		\enspace.
		\end{align*}
		Summing over $t$, we conclude the proof.
	\end{proof}
	
			\subsection{Exponential moments}\label{sec:expmom}
	
	The following proposition is the central piece of the proof of \Cref{thm:main}.
	
	\begin{proposition}\label{prop:expmom}
		Assume \ref{ass:normalized}, \ref{ass:erg}, \ref{ass:rate} and \ref{ass:domination}. There is $\theta_\star>0$ such that, for every $\theta\in(0,\theta_\star]$, there are $\varrho_\theta\in(0,1)$ and $C_\theta<\infty$ with, for every $z\in\Zcal$, every $t\ge0$ and every function $\varphi:\Zcal\to\R_{++}$ such that $\norm{\varphi}_V<+\infty$:
		\begin{equation}\label{eq:expmom}
			\E_z\big[e^{\theta S^0_t}\varphi(z^0_t)\big]\le C_\theta\varrho_\theta^{t}\norm{\varphi}_V V(z) \enspace,
		\end{equation}
		where $\E_z$ refers to the limit chain $\{z_t^0\}$ initialized at $z$
		and $S_t^0 = \sum_{k<t} \log\Gamma_0(z_k^0,\mathbf U_{k+1})$.
	\end{proposition}
	\begin{proof}
		This proof uses a perturbation argument for quasi-compact operators on a Banach space, from \cite{hennion2001limit}.
		Throughout the proof, $\mathcal B:=\mathsf L^V$ is the complex Banach space of measurable $\varphi:\Zcal\to\mathbb C$ with $\norm\varphi_V:=\sup\abs\varphi/V<\infty$, and $\mathcal L_{\mathcal B}$ is the space of bounded linear operators on $\mathcal B$, equipped with the operator norm $\norm{\cdot}_{V\to V}$.
		For $\theta \in \mathsf B(0,\Lambda) = \{z\in\mathbb C : \abs z \le \Lambda\}$, where $\Lambda$ is given in \ref{ass:domination}, we define the operator $\widehat P_\theta$ on $\mathcal B$ by
		\[
		\widehat P_\theta\varphi(z):=\int\Gamma_0(z,u)^{\theta}\varphi\big(\Phi_0(z,u)\big)\nu(\dd u) \enspace,
		\]
		so that $\widehat P_0=P_0$ and, by the Markov property, $\widehat P_\theta^{t}\varphi(z)=\E_z\big[e^{\theta S^0_t}\varphi(z^0_t)\big]$ for every $t\ge0$.
		
		In order to apply the perturbation theorem \cite[Theorem III.8]{hennion2001limit}, we first check its two hypotheses:
		(H1) $\theta\mapsto\widehat P_\theta$ is  $\mathcal C^2$ from $\mathsf B(0,\Lambda)$ to $\mathcal L_{\mathcal B}$, and
		(H2) $P_0$ has one dominating simple eigenvalue (see \cite[Definition III.7]{hennion2001limit}) and its spectral radius $r(P_0)$ equals one.

		\emph{(i) Proof of (H1).} For $k\ge0$ let $\widehat P^{(k)}\varphi(z):=\int\big(\log\Gamma_0(z,u)\big)^k\varphi(\Phi_0(z,u))\nu(\dd u)$. By \ref{ass:normalized} and \ref{ass:domination}, $\abs{\log\Gamma_0}^k\le\ell^k$ and $V(\Phi_0(z,u))\le V(z)\Xi(u)$, so that $\norm{\widehat P^{(k)}}_{V\to V}\le\int\ell^k\Xi\dd\nu$, whence, for $\abs\theta\le \Lambda$,
		\[
		\sum_{k\ge0}\frac{\abs\theta^k}{k!}\norm{\widehat P^{(k)}}_{V\to V}\le\int e^{\abs\theta\ell}\Xi\dd\nu<\infty \enspace,
		\]
		by \ref{ass:domination}. 
		The series $\widehat P_\theta=\sum_{k\ge0}\frac{\theta^k}{k!}\widehat P^{(k)}$ therefore converges in operator norm and defines a $\mathcal C^\infty$ family of bounded operators on $\mathsf L^V$ on  $\mathsf B(0,\Lambda)$.
		Note that in particular, its derivative in $\theta=0$ is $\widehat P^{(1)}\ind=\gamma_0$ by \eqref{eq:def-P0}.

		\emph{(ii) Proof of (H2).}
		Let $\Pi:=\ind\otimes\pi_0$, that is $\Pi\varphi=\pi_0(\varphi)\ind$, where $\pi_0$ is the invariant probability measure of $P_0$ given by \ref{ass:erg}. It is a bounded operator on $\mathcal B$, since $\abs{\pi_0(\varphi)}\le\pi_0(V)\norm\varphi_V$ with $\pi_0(V)<\infty$, and a projection, since $\pi_0(\ind)=1$. We check the four conditions (i)--(iv) of \cite[Definition III.7]{hennion2001limit} for the closed subspaces $\mathsf F:=\mathrm{range}\,\Pi$ and $\mathsf G:=\ker\Pi$, and we prove that $r(P_0)=1$.
		As $\Pi$ is a bounded projection, $\mathcal B=\mathsf F\oplus\mathsf G$, which is condition (i).
		As $P_0\ind=\ind$, the subspace $\mathsf F$ is invariant under $P_0$, and as $\pi_0P_0=\pi_0$, every $\varphi\in\mathsf G$ satisfies $\pi_0(P_0\varphi)=\pi_0(\varphi)=0$, so that $\mathsf G$ is invariant as well, which is condition (ii).
		The subspace $\mathsf F=\mathbb C\ind$ is one-dimensional and $P_0$ acts on it as the identity, so that $P_0\vert_{\mathsf F}$ has the single simple eigenvalue $1$, which is condition (iii).
		Finally, $P_0^n\varphi=(P_0^n-\Pi)\varphi$ for every $\varphi\in\mathsf G$, so that $\norm{P_0^n\vert_{\mathsf G}}_{V\to V}\le M\rho_V^n$ by \eqref{eq:Vgeo} and $r(P_0\vert_{\mathsf G})\le\rho_V<1$, which is condition (iv).
		The same bound gives $\sup_n\norm{P_0^n}_{V\to V}<\infty$, so $r(P_0)\le1$, and $r(P_0)\ge1$ since $P_0\ind=\ind$.
		Therefore (H2) is verified.

		Apply the perturbation theorem \cite[Theorem III.8]{hennion2001limit} with $m=2$ and one dominating eigenvalue. 
		There are therefore $\theta_1\in(0,\Lambda)$, $\varrho_-\in[0,1)$ and $c<\infty$ such that, for $\abs\theta<\theta_1$, the operator $\widehat P_\theta$ still has one dominating simple eigenvalue $\lambda_\theta$, of class $\mathsf C^2$ in $\theta$, with $\lambda_0=1$ and $\varrho_\theta:=\abs{\lambda_\theta}>\varrho_-$, with right and left eigenvectors $v_\theta\in\mathsf L^V$ and $\mu_\theta$ normalized by $\mu_\theta(v_\theta)=1$, and with a remainder $N_\theta$ such that
		\begin{equation}\label{eq:HH}
			\widehat P_\theta^{t}=\lambda_\theta^{t}v_\theta\otimes\mu_\theta+N_\theta^{t},\qquad \norm{N_\theta^{t}}_{V\to V}\le c\varrho_-^{t}\qquad(t\ge1) \enspace.
		\end{equation}
		The derivative of $\lambda_\theta$ at $0$ is
		\[
		\frac{\dd}{\dd\theta}\lambda_\theta\Big\vert_{\theta=0}=\pi_0\big(\widehat P^{(1)}\ind\big)=\pi_0(\gamma_0)=-\CR<0
		\]
		by \cite[Corollary III.11]{hennion2001limit} and \ref{ass:rate}. Hence $\lambda_\theta=1-\theta\CR+o(\theta)$ and $\varrho_\theta\le\abs{1-\theta\CR}+o(\theta)<1$ for $\theta>0$ small, while $\varrho_\theta>\varrho_-\ge0$. There is therefore $\theta_\star\in(0,\theta_1)$ such that $0<\varrho_\theta<1$ for every $\theta\in(0,\theta_\star]$.
		
		Let $\theta\in(0,\theta_\star]$. Since $\abs{\lambda_\theta}=\varrho_\theta>\varrho_-$, \eqref{eq:HH} gives $\norm{\widehat P_\theta^{t}}_{V\to V}\le C_\theta\varrho_\theta^{t}$ for every $t\ge0$, with $C_\theta:=\max\big(1,\norm{v_\theta}_V\norm{\mu_\theta}+c\big)$, the case $t=0$ being trivial.
		This ends the proof.
	\end{proof}

	\begin{remark}
		What has been checked in the proof of (H2) is that $P_0$ is $V$-uniform in the sense of \cite[Section 4.1]{kontoyiannis2003spectral}, see \cite[Proposition 4.6]{kontoyiannis2003spectral}, whose condition (iii) is \eqref{eq:Vgeo}. It is a consequence of \ref{ass:erg} and not a property of Markov kernels in general: it fails as soon as $P_0$ is reducible, periodic, or ergodic at a subgeometric rate.
	\end{remark}

			\subsection{Probability of the identification event}\label{sec:ident}
	
	The exponential moments in \Cref{prop:expmom} allow us to bound the probability given by \Cref{prop:bound-Ec}. This is summarized in the next proposition.
	
	\begin{proposition}\label{prop:ident}
		Assume \ref{ass:normalized}--\ref{ass:domination}.
		Let $\theta:=\min(\theta_\star,\alpha,\eta)$, with $\theta_\star$ as in \Cref{prop:expmom} and $\alpha$, $\eta$ as in \ref{ass:coupling}, let $\varrho_\theta,C_\theta$ be the constants of \Cref{prop:expmom} for this $\theta$, and let
		\begin{equation}\label{eq:C}
			C:=\frac{C_\theta}{1-\varrho_\theta}\Big(\kappa\bar\sigma^{\eta-\theta}+\bar\sigma^{-\theta}+1\Big) \enspace.
		\end{equation}
		Then
		the event $\Ecal$ in \eqref{eq:Eevent} satisfies $\Prob(\Ecal^c)\le CV(z_0)\sigma_0^{\theta}$.
	\end{proposition}
	\begin{proof}
		We bound the three sums in \eqref{eq:bound-E}. 
		
		Since $\sigma^0_t=\sigma_0e^{S^0_t}$ and $\ind\le V$, \Cref{prop:expmom} gives, for every $t\ge0$,
		\begin{equation}\label{eq:twomoments}
		\begin{aligned}
			\E\big[(\sigma^0_t)^{\theta}\big]&\le\sigma_0^{\theta}\E\big[e^{\theta S^0_t}\big]\le C_\theta\varrho_\theta^{t}V(z_0)\sigma_0^{\theta} \enspace,\\
			 \E\big[V(z^0_t)(\sigma^0_t)^{\theta}\big]&=\sigma_0^{\theta}\E\big[e^{\theta S^0_t}V(z^0_t)\big]\le C_\theta\varrho_\theta^{t}V(z_0)\sigma_0^{\theta} \enspace.
		\end{aligned}
		\end{equation}
		As $\theta\le\eta$, we have $(\sigma_t^0)^{\eta}\ind_{\{\sigma_t^0\le\bar\sigma\}}\le\bar\sigma^{\eta-\theta}(\sigma_t^0)^{\theta}$. 
		Hence, we have
		 $$
		 \kappa \sum_{t=0}^{+\infty} \E[
		 \ind_{\{\sigma^0_t\le\bar\sigma\}}(\sigma_t^0)^\eta]
		 \le
		 \kappa\bar\sigma^{\eta-\theta}\sum_{t\ge0}\E[(\sigma^0_t)^{\theta}]
		 \le 
		 \kappa \bar\sigma^{\eta-\theta} C_\theta V(z_0) \sigma_0^\theta \times (1-\varrho_\theta)^{-1}
		  \enspace.
		 $$
		 which bounds the first sum in \eqref{eq:bound-E}.

		 As $\theta\le\alpha$ and $V\ge1$, we have $V^{\theta/\alpha}\le V$. 
		 Hence, when $V(z^0_t)>(\sigma^0_t)^{-\alpha}$ is true, it implies $V(z^0_t)(\sigma^0_t)^{\theta}>1$. 
		 Markov's inequality gives then
		 $
		 \Prob\Big(V(z^0_t)>\big(\sigma^0_t\big)^{-\alpha}\Big)\le\E\big[V(z^0_t)(\sigma^0_t)^{\theta}\big] 
		 $.
		 The second sum in \eqref{eq:bound-E} therefore satisfies 
		 $$
		  \sum_{t=0}^{+\infty}\Prob [V(z_t^0)>(\sigma_t^0)^{-\alpha}] \le
		 \sum_{t\ge0}\E[V(z^0_t)(\sigma^0_t)^{\theta}]
		 \le  C_\theta V(z_0)\sigma_0^{\theta} \times (1-\varrho_\theta)^{-1} 
		  \enspace.
		 $$
		 For the third sum, 
		 by Markov's inequality,
		$
		\Prob\big(\sigma^0_t>\bar\sigma\big)\le\bar\sigma^{-\theta}\E\big[(\sigma^0_t)^{\theta}\big] 
		$,
		so that 
		$$
		\sum_{t=0}^{+\infty} \Prob[ \sigma_t^0>\bar\sigma ] \le
		\bar\sigma^{-\theta}\sum_{t\ge0}\E[(\sigma^0_t)^{\theta}]
		\le \bar\sigma^{-\theta} \times C_\theta V(z_0)\sigma_0^{\theta} \times (1-\varrho_\theta)^{-1}
		 \enspace.
		$$ 
		Summing the three bounds, we obtain the desired bound, with $C$ as in \eqref{eq:C}.
	\end{proof}
	
	\subsection{Convergence of the limit chain and proof of \Cref{thm:main}}\label{sec:proofmain-conclusion}
	
	Before concluding the proof, we give in the next lemma limit results on the limit chain \eqref{eq:limit-chain}.
	
	\begin{lemma}\label{lem:limit-asymp}
		Assume \ref{ass:normalized}--\ref{ass:domination}.
		Almost surely, $S^0_t/t\to-\CR$ and $\frac1t\log V(z^0_t)\to0$. 
	\end{lemma}
	\begin{proof}
		By \ref{ass:domination}, the function $e^{\Lambda\ell}$ is $\nu$-integrable, since $\Xi\ge1$.
		In particular $\ell\in\mathsf L^2(\nu)$ and, by \ref{ass:normalized}, $\abs{\gamma_0}\le\int\ell\dd\nu<\infty$ on $\Zcal$.
		Define
		\[
		\chi_{r+1}=\log\Gamma_0(z^0_r,\mathbf U_{r+1})-\gamma_0(z^0_r) \enspace.
		\]
		We have $\E[\chi_{r+1}\mid\Fcal_r]=0$ by \eqref{eq:def-P0}. 
		Hence $\sum_{r<t}\chi_{r+1}$ is an $(\Fcal_t)$-martingale.
		Besides, by \ref{ass:normalized}, $\abs{\chi_{r+1}}\le\ell(\mathbf U_{r+1})+\int\ell\dd\nu$.
		Since $\ell\in\mathsf L^2(\nu)$, $\E[\chi_{r+1}^2]\le\int(\ell+\int\ell\dd\nu)^2\dd\nu<\infty$ for every $r$, so $\sum_{r\ge0}\E[\chi_{r+1}^2]/(r+1)^2<\infty$.
		By the strong law of large numbers for martingales \cite{chow1967strong}, $\frac1t\sum_{r<t}\chi_{r+1}\to0$ almost surely.
		
		 Moreover, by \ref{ass:erg} the chain $\{z^0_t\}_{t\in\N}$ is geometrically ergodic with invariant probability $\pi_0$~\cite[Theorem  15.0.1]{meyn2009markov}, and it is positive Harris recurrent: by \cite[Lemma 15.2.2]{meyn2009markov} the drift of \ref{ass:erg} implies that $V$ is unbounded off petite sets, so the chain is Harris recurrent by \cite[Theorem 9.1.8]{meyn2009markov}, and it has an invariant probability.
		 By the law of large numbers for positive Harris chains~\cite[Theorem 17.1.7]{meyn2009markov}, since $\gamma_0$ is bounded, we have almost surely
		 $$\lim_{t\to\infty}\frac1t\sum_{r<t}\gamma_0(z^0_r)  =  \pi_0(\gamma_0)=-\CR \enspace.$$
		 Summing the two limits, $\lim S^0_t/t=-\CR$.
		
		Since $V\ge1$, $\liminf_t\frac1t\log V(z^0_t)\ge0$. 
		Besides, iterating \ref{ass:erg} gives 
		$$\E V(z^0_t)\le\lambda_V^tV(z_0)+b_V/(1-\lambda_V)\le V(z_0)+b_V/(1-\lambda_V) \enspace.$$
		Hence, for every $\varepsilon>0$, 
		$$\Prob(V(z^0_t)>e^{\varepsilon t})\le\big(V(z_0)+b_V/(1-\lambda_V)\big)e^{-\varepsilon t}$$
		 is summable for $t\in\N$.
		By Borel--Cantelli, $\limsup_t\frac1t\log V(z^0_t)\le\varepsilon$ almost surely. 
	\end{proof}
	
	We can now conclude the proof of \Cref{thm:main}.
	
	\begin{proof}[Proof of \Cref{thm:main}]
		Let $\theta$ and $C$ be as in \Cref{prop:ident}, and let $(z_0,\sigma_0)\in\Zcal\times(0,+\infty)$. Let $\Ecal$ be the event of \eqref{eq:Eevent}.
		Then $\Prob(\Ecal^c)\le CV(z_0)\sigma_0^{\theta}$ by \Cref{prop:ident}, and on $\Ecal$ we have $(z_t,\sigma_t)=(z^0_t,\sigma^0_t)$ for all $t\ge0$. 
		Hence $\frac1t\log\sigma_t=\frac1t\log\sigma_0+\frac1tS^0_t$ and $V(z_t)=V(z^0_t)$ on $\Ecal$, so \Cref{lem:limit-asymp} gives the two limits there, up to a null set. 
		
		For the marginal law, we use the convention $\norm{\mu-\mu'}_\TV:=\sup_{\abs\varphi\le1}\abs{\mu(\varphi)-\mu'(\varphi)}$ of \cite{meyn2009markov}, which is twice the supremum over measurable sets. Let $\abs\varphi\le1$. Then $\abs{\E\varphi(z_t)-\E\varphi(z^0_t)}\le2\Prob(z_t\ne z^0_t)\le2\Prob(\Ecal^c)$, and $\Lcal(z^0_t)=P_0^t(z_0,\cdot)$, which ends the proof.
		
	\end{proof}

	\section{Proof of \Cref{th:decoupled-initial-conditions}}\label{sec:proofdecoupled}

	We prove in this section \Cref{th:decoupled-initial-conditions}.
	The following lemma uses results of the drift theory of \cite{hajek1982hitting} to obtain an occupation bound useful to quantify the time which the adaptive chain \eqref{eq:adaptive-z} needs to reach a convergence area covered by \Cref{thm:main}.
	
	\begin{lemma}\label{l:occupation}
		Assume \ref{assB:drift}.
		Then $\tau<\infty$ almost surely and
		\begin{equation}\label{e:occ}
			\E\sum_{t<\tau}e^{-\vartheta D(z_t)}\le C_{\rm occ}:=\frac{e^{\delta}}{(1-\rho)(1-e^{-\vartheta/2})} \enspace,
		\end{equation}
		where $\delta\in(0,\vartheta/2]$ and $\rho\in(0,1)$ are constants which depend only on $(c,\vartheta,K)$.
	\end{lemma}
	
	\begin{proof}
		The inequalities \eqref{eq:drift-for-zsig} are true for $t<\tau$ only. We set
		\begin{equation*}
			\hat\Delta_{t+1}:=(D(z_t)-D(z_{t+1}))\ind_{\{t<\tau\}}+c\ind_{\{t\ge\tau\}},
			\qquad \hat D_t:=D(z_0)-\sum_{u=1}^{t}\hat\Delta_u \enspace.
		\end{equation*}
		Then $\hat D_t=D(z_t)$ for $t\le\tau$. 
		Besides, we have for all $t\ge0$:
		\begin{equation}\label{eq:inproof-occupation-1}
			\E[\hat \Delta_{t+1} |\Fcal_t] \ge c , \qquad \E[e^{\vartheta|\hat\Delta_{t+1}|} | \Fcal_t ]  \le K .  
		\end{equation}
		 We check the two conditions C1 and C2 of the drift analysis of \cite{hajek1982hitting} for the sequence $(\hat D_t,\mathcal F_t)_{t\ge0}$ and the level $0$. 
		 Condition C1 holds with $\varepsilon_0=c$ by \eqref{eq:inproof-occupation-1}. 
		 
		 By Markov's inequality, $\Prob(\abs{\hat\Delta_{t+1}}>x\mid\mathcal F_t)\le\min(1,Ke^{-\vartheta x})$ for all $x\ge0$. 
		 Define $Z$ a random variable on $\R_{+}$ by $\Prob(Z>x)=\min(1,Ke^{-\vartheta x})$ for $x\ge 0$.
		 Then, 
		 \begin{multline*}
		 	\E e^{\vartheta Z/2} = \int_{y\in\R_+} \Prob[e^{\vartheta Z/2}>y] \dd y = \int_{y\in\R_+} \min\{1, K y^{-2} \} \dd y \\
		 	= K^{1/2} + \int_{K^{1/2}}^{+\infty} K y^{-2} \dd y = 2 K^{1/2} < +\infty.
		 \end{multline*}
		 Hence Condition C2 holds with majorant $Z$. Lemma 2.1 of \cite{hajek1982hitting} then provides $\delta\in(0,\vartheta/2]$ and $\rho\in(0,1)$, depending only on $(c,\vartheta,K)$, such that Conditions D1 and D2 hold with the exponent $\delta$. 
		 Let $\hat\tau_0:=\inf\{t\ge0:\hat D_t\le0\}$. The hitting-time bound (2.9) of \cite[Theorem 2.3]{hajek1982hitting} reads $\Prob(\hat\tau_0>k\mid\mathcal F_0)\le e^{\delta\hat D_0}\rho^k$. 
		 Since $\hat D_k=D(z_k)>0$ on $\{\tau>k\}$, we have $\{\tau>k\}\subset\{\hat\tau_0>k\}$. The conditions hold at every time, so the same bound applies to the sequence $(\hat D_{t+k},\mathcal F_{t+k})_{k\ge0}$, and for all $t\le k$
		\begin{equation}\label{e:hit}
			\Prob(\tau>k\mid\mathcal F_t)\le\rho^{k-t}e^{\delta\hat D_t} \enspace,
		\end{equation}
		In particular, $\Prob(\tau>k)\le\rho^ke^{\delta D(z_0)}\to0$, so $\tau<\infty$ almost surely which proves the first claim of \Cref{l:occupation}.
		
		We prove now \eqref{e:occ}.
		For $t<\tau$ we have $D(z_t)>0$. 
		Then,
		\begin{equation*}
			\E\sum_{t<\tau}e^{-\vartheta D(z_t)}\le\sum_{j\ge0}e^{-\vartheta j}\E\#\{t<\tau \mid j\le D(z_t)<j+1\}  \enspace.
		\end{equation*}
		However
		\begin{equation*}
			\E\#\{t<\tau \mid  D(z_t)<j+1\} 
			\le \E [\ind_{\hat\tau_{j+1}<\tau} (\tau-\hat\tau_{j+1})]
			= \sum_{t=0}^{+\infty} \E\left[   \ind_{t=\hat\tau_{j+1}< \tau} (\tau-t)   \right]
		\end{equation*}
		where $\hat\tau_{j+1}=\inf\{  t \mid \hat D_t \le j+1 \}$.
		Yet, summing \eqref{e:hit} for $k= t,t+1,\dots$ gives $\E[\ind_{t<\tau}(\tau-t)\mid\mathcal F_t]\le e^{\delta\hat D_t}/(1-\rho)$.
		Hence,
		\begin{equation*}
			\E\#\{t<\tau \mid  D(z_t)<j+1\} 
			\le
			\sum_{t=0}^{+\infty} \E\left[   \ind_{t=\hat\tau_{j+1}< \tau} \frac{ e^{\delta\hat D_t}}{1-\rho}  \right]
			\le \frac{ e^{\delta(j+1)}}{1-\rho} \enspace.
		\end{equation*}
		Summing over $j$:
		\begin{equation*}
			\E\sum_{t<\tau}e^{-\vartheta D(z_t)}\le \sum_{j\ge0} \frac{ e^{\delta(j+1)-\vartheta j}}{1-\rho}\le\frac{e^{\delta}}{(1-\rho)(1-e^{-\vartheta/2})}
		\end{equation*}
		since $\delta\le \vartheta/2$. This ends the proof.
	\end{proof}

	We have now the required tools to prove \Cref{th:decoupled-initial-conditions}.

		\begin{proof}[Proof of \Cref{th:decoupled-initial-conditions}]
			By \Cref{l:occupation}, $\tau$ is an almost surely finite stopping time, so the samples $(\mathbf U_{\tau+t})_{t\ge1}$ are i.i.d.\ with law $\nu$ and independent of $\Fcal_\tau$.
			Hence \Cref{thm:main}, applied to the process $(z_{\tau+t},\sigma_{\tau+t})_{t\ge0}$ started from $(z_\tau,\sigma_\tau)$, provides an event $\Ecal_\tau$ with
			\begin{equation}
				\Prob[\Ecal_\tau^c | \Fcal_\tau] \le \min\{ 1, C V(z_\tau) \sigma_\tau^\theta  \}  \enspace,
			\end{equation}
			where $C,\theta$ are positive constants given by \Cref{thm:main}. Without loss of generality, we assume the constant $\vartheta$ of \ref{assB:drift} satisfies $\vartheta\le\theta$.
			Indeed, if $\vartheta>\theta$, Jensen's inequality gives $\E[e^{\theta\abs{D(z_t)-D(z_{t+1})}}\mid\Fcal_t]\le K^{\theta/\vartheta}$, and $K^{\theta/\vartheta}\ge e^{\theta c}$, so \ref{assB:drift} holds with $(\theta,K^{\theta/\vartheta})$ in place of $(\vartheta,K)$, and the constant $C_{\rm occ}$ of \Cref{l:occupation} is then computed with these values. This is why $C_{\rm occ}$ depends on $\theta$ in the statement.
			 We have that
			\begin{equation}\label{eq:inproof-th2-1}
				\E[\sigma_\tau^\theta] \le \sigma_0^\theta + \E\left[ \sum_{t<\tau} \sigma_{t+1}^\theta \right] \enspace.
			\end{equation}
			However, using \ref{ass:normalized} and \ref{assB:drift}, we have for $t < \tau$
			\begin{equation}\label{eq:inproof-th2-2}
				\E[\sigma_{t+1}^\theta | \Fcal_t]  
				= \sigma_t^\theta \times \int \Gamma_{\sigma_t}(z_t,u)^\theta \nu(\dd u) 
				\le \varsigma^\theta e^{-\theta D(z_t)} \int e^{\theta \ell} \dd \nu
				\le  \varsigma^\theta e^{-\vartheta D(z_t)} \int e^{\theta \ell} \dd \nu \enspace,
			\end{equation}
			where we used for the last inequality that $\vartheta\le\theta$ and $D(z_t)\ge 0$ for $t<\tau$.
			The integral $\int e^{\theta\ell}\dd\nu$ is finite by \ref{ass:domination}, since $\Xi\ge1$ and $\theta\le\theta_\star<\theta_1<\Lambda$ by \Cref{prop:ident} and the proof of \Cref{prop:expmom}.
			Combining \eqref{eq:inproof-th2-1} and \eqref{eq:inproof-th2-2}, we get,  using the constant $C_{\rm occ}$ from \Cref{l:occupation}:
			\begin{equation}
				\E[\sigma_\tau^\theta] \le \sigma_0^\theta + C_{\rm occ} \varsigma^\theta \int e^{\theta \ell} \dd\nu \enspace.
			\end{equation}
			By definition of $\tau$ in \ref{assB:drift}, we have either $D(z_\tau)\le0$ or $g(z_\tau,\sigma_\tau)>1$.
			On the one hand, on $(D(z_\tau)\le0)$, we have $V(z_\tau)\le \sup_{D\le0} V$.
			On the other hand, $\Prob[g(z_\tau,\sigma_\tau)>1]\le \E[g(z_\tau,\sigma_\tau)]\le g(z_0,\sigma_0)$ by the third inequality in \ref{assB:drift}.
			Indeed, this inequality says that the stopped process $\{g(z_{t\wedge\tau},\sigma_{t\wedge\tau})\}_{t\in\N}$ is a nonnegative supermartingale, so $\E[g(z_{t\wedge\tau},\sigma_{t\wedge\tau})]\le g(z_0,\sigma_0)$ for every $t$, and Fatou's lemma with $\tau<\infty$ almost surely gives $\E[g(z_\tau,\sigma_\tau)]\le g(z_0,\sigma_0)$.
			All in all, we get:
			\begin{multline*}
				\Prob[\Ecal_\tau^c] \le
				\Prob[g(z_\tau,\sigma_\tau)>1]
				+ \E [  C V(z_\tau) \sigma_\tau^\theta  \ind_{ D(z_\tau)\le0 } ] 
				\le g(z_0,\sigma_0) + C \sup_{D\le0} V \times  \left(  \sigma_0^\theta + C_{\rm occ} \varsigma^\theta \int e^{\theta \ell} \dd\nu  \right)
			\end{multline*}
			which ends the proof.
		\end{proof}
	
	\section{Proof of \Cref{thm:es2}}\label{sec:proof-es}

		\subsection{Notations}\label{sec:es-notation}

		For $u=(u^1,\dots,u^\lambda)\in\Ucal$ and  $s\in\Sym_\lambda$ a permutation of $\{1,\dots,\lambda\}$, let $Y_s(u):=\sum_{i\le\mu}w_iu^{s(i)}$.
		We write $\abs u_\infty:=\max_{i\le\lambda}\abs{u^i}$, so that $\abs{Y_s(u)}\le\abs u_\infty$ since the weights sum to one.
		For  $F:\R^d\to\R$, let $s_F(u)\in\Sym_\lambda$ be the ranking of $u^1,\dots,u^\lambda$ by $F$, ties being broken by the index of the candidate, so that $s_F(u)$ is uniquely defined for every $u$.

		The normalized process is $z_t:=(m_t-x^*)/\sigma_t\in\Zcal:=\R^d$. 
		Since $f(m_t+\sigma_tu)=f\big(x^*+\sigma_t(z_t+u)\big)$, the selected step of \eqref{eq:es} is $Y_{t+1}=Y_{f_{\sigma_t}(z_t+\cdot)}(\mathbf U_{t+1})$ with $f_\sigma(x):=f(x^*+\sigma x)$ for $\sigma>0$, and \eqref{eq:es} takes the form \eqref{eq:normalized-dynamics-sigma} with
		\begin{equation}\label{eq:esmaps}
			\Phi_\sigma(z,u):=\frac{z+c_mY_{f_\sigma(z+\cdot)}(u)}{\Gamma\big(Y_{f_\sigma(z+\cdot)}(u)\big)},\qquad \Gamma_\sigma(z,u):=\Gamma\big(Y_{f_\sigma(z+\cdot)}(u)\big) \enspace,
		\end{equation}
		where $\Gamma$ is the step-size change \eqref{eq:CSA2}.
		The limit map and the limit rule are those of the algorithm on the quadratic model $f_0$, namely
		\begin{equation}\label{eq:limitmaps}
			\Phi_0(z,u):=\frac{z+c_mY_{f_0(x^*+z+\cdot)}(u)}{\Gamma\big(Y_{f_0(x^*+z+\cdot)}(u)\big)},\qquad \Gamma_0(z,u):=\Gamma\big(Y_{f_0(x^*+z+\cdot)}(u)\big) \enspace.
		\end{equation}
		Here $f_0$ is the quadratic model of $f$ and not the function $f_\sigma$ at $\sigma=0$, which is constant.
	
		\subsection{Verification of \ref{ass:normalized}--\ref{ass:domination}}
		
		\begin{proposition}\label{prop:esA1}
			The multiplicative rule \eqref{eq:normalized-dynamics-sigma} of \ref{ass:normalized} holds.
		\end{proposition}
		\begin{proof}
			Straightforward.
		\end{proof}
		
		\begin{proposition}\label{prop:esA2}
			For every $\sigma>0$, the maps \eqref{eq:esmaps} of the algorithm \eqref{eq:es} on $f_0(x)=(x-x^*)^\top \nabla^2 f(x^*)(x-x^*)$ are the limit maps \eqref{eq:limitmaps}. In particular, $P_0$ is the transition kernel of the normalized process of the algorithm on $f_0$.
		\end{proposition}
		\begin{proof}
			The function $f_0$ is scaling-invariant. By \cite[Proposition 2]{toure2023global}, the ranking of the candidates by $f_0$ does not depend on $\sigma$, so the ranking of $z+u^1,\dots,z+u^\lambda$ by $f_0(x^*+\sigma\cdot)$ is the ranking by $f_0(x^*+\cdot)$, and the maps \eqref{eq:esmaps} for $f_0$ are \eqref{eq:limitmaps}.
		\end{proof}
		
		\begin{proposition}\label{prop:esA3}
			Let $\beta\in(0,\min(1,\beta_0)]$ as in \eqref{eq:Vbeta}.
			The domination of \ref{ass:erg} holds with
			\begin{equation*}
				V_\beta(z)=1+\abs z^\beta\qquad\text{and}\qquad\Xi_\beta(u)=\max_{s\in\mathfrak{S}_\lambda}\big[1+ \Gamma(Y_s(u))^{-\beta} +\Gamma(Y_s(u))^{-\beta} c_m^\beta \abs{Y_s(u)}^\beta\big] \enspace.
			\end{equation*}
		\end{proposition}
		\begin{proof}
			We have 
			\begin{align*}
				V_{\beta} (\Phi_\sigma(z,u)) 
				& = 1 + \left| \frac{z+c_m Y_{s_{z,\sigma}}(u)}{\Gamma( Y_{s_{z,\sigma}}(u))} \right|^\beta \\
				& \le (1+ \abs z^\beta) \times \left[ 1+ \Gamma( Y_{s_{z,\sigma}}(u))^{-\beta} + \Gamma( Y_{s_{z,\sigma}}(u))^{-\beta} c_m^\beta\abs{   Y_{s_{z,\sigma}}(u) }^\beta  \right]  \\
				& \le V_\beta(z) \times \Xi_\beta(u) \enspace.
			\end{align*}
		\end{proof}
		
		\begin{proposition}\label{prop:esA4}
			The bound $\abs{\log\Gamma_\sigma(z,u)}\le\ell(u)$ of \ref{ass:normalized} is true with $\ell(u)=\frac{1}{2d_\sigma}\Big(1+\frac{\abs u_\infty^2}{d\abs w^2}\Big)$.
		\end{proposition}
		\begin{proof}
			By \eqref{eq:CSA2} and \eqref{eq:esmaps}, $\log\Gamma_\sigma(z,u)=\frac{1}{2d_\sigma}\big(\frac{\abs{Y_s(u)}^2}{d\abs w^2}-1\big)$ for some $s\in\Sym_\lambda$, and $\abs{Y_s(u)}\le\abs u_\infty$.
		\end{proof}

		\begin{proposition}\label{prop:esA6}
			Assume \ref{ass:function} and let $\beta\in(0,\min(1,\beta_0)]$ as in \eqref{eq:Vbeta}. Then \ref{ass:coupling} holds with $V=V_\beta$, $\alpha=\beta/4$ and $\eta=1/8$.
		\end{proposition}
		\begin{proof}
			Let $L$ be a Lipschitz constant of $\nabla^2f$ on $B(x^*,r_0)$, with $r_0$ as in \ref{ass:function}, up to taking $r_0$ smaller if necessary.
			Let $q(x)=x^\top\nabla^2f(x^*)x$, so that $f_0(x^*+\sigma x)=\sigma^2q(x)$.
			Taylor's formula  gives, for $\sigma\abs x<r_0$,
			\begin{equation}\label{eq:taylor}
				f_\sigma(x)=f(x^*)+\frac{\sigma^2}{2}\big(q(x)+\omega_\sigma(x)\big)\qquad\text{with}\qquad\abs{\omega_\sigma(x)}\le\frac L3\sigma\abs x^3 \enspace.
			\end{equation}
			Let $\bar\sigma\in(0,1)$ with $2\bar\sigma^{3/4}<r_0$, let $\sigma\le\bar\sigma$ and $z\in\Zcal_\sigma$. Then $\abs z^\beta\le V_\beta(z)\le\sigma^{-\beta/4}$, that is $\abs z\le\sigma^{-1/4}$.

			By \eqref{eq:esmaps}, $\Phi_\sigma$ and $\Gamma_\sigma$ are functions of the ranking $s_\sigma$ of $z+u^1,\dots,z+u^\lambda$ by $f_\sigma$. 
			Moreover, $\Phi_0$ and $\Gamma_0$ are the same functions of the ranking $s_0$ by $q$. 
			Hence $u\in\Dcal_\sigma(z)$ implies that $s_\sigma\neq s_0$. 
			If $\max_i\abs{u^i}\le\sigma^{-1/4}$, then $\sigma\abs{z+u^i}\le2\sigma^{3/4}<r_0$ for all $i$, so \eqref{eq:taylor} applies for $x=z+u^i$. 
			The ranking $s_\sigma$ by $f_\sigma$ is then the ranking by $q+\omega_\sigma$. 
			Yet $s_0\neq s_\sigma$ when some pair $i\ne j$ is ordered differently, i.e., $q(z+u^i)\le q(z+u^j)$ and $q(z+u^j)+\omega_\sigma(z+u^j)\le q(z+u^i)+\omega_\sigma(z+u^i)$. 
			Therefore
			\begin{equation*}
				\abs{q(z+u^i)-q(z+u^j)}\le\abs{\omega_\sigma(z+u^i)}+\abs{\omega_\sigma(z+u^j)}\le\frac{2L}{3}\sigma\big(2\sigma^{-1/4}\big)^3=\frac{16L}{3}\sigma^{1/4} \enspace.
			\end{equation*}
			We have shown
			\begin{equation}\label{eq:Dsplit}
				\nu\big(\Dcal_\sigma(z)\big)\le\nu\Big(\max_i\abs{u^i}>\sigma^{-1/4}\Big)+\sum_{i<j}\nu\Big(\abs{q(z+u^i)-q(z+u^j)}\le\frac{16L}{3}\sigma^{1/4}\Big) \enspace.
			\end{equation}

			For the first term, $\int e^{\abs{u^i}^2/4}\dd\nu=2^{d/2}$ and Markov's inequality give 
			$$
			\nu(\abs{u^i}>\sigma^{-1/4}) = \nu\left( e^{\abs{u^i}^2/4} > e^{\sigma^{-1/2}/4} \right)
			\le2^{d/2}e^{-\sigma^{-1/2}/4}\le2^{d/2}\sigma^{1/8}\enspace.
			$$
			Summing over $i$:
			\begin{equation*}
				\nu\Big(\max_i\abs{u^i}>\sigma^{-1/4}\Big)\le\lambda2^{d/2}\sigma^{1/8} \enspace.
			\end{equation*}

			For the second term we use the Carbery--Wright inequality \cite[Theorem 8]{carbery2001distributional}, applied to the Gaussian measure with a polynomial of degree at most $2$ and with the exponent $q=4$.
			It gives a constant $C_{\rm CW}>0$, independent of $n$, such that, for any quadratic function $Q$ on $\R^n$, every $X\sim\mathcal N(0,I_n)$ and every $\varepsilon>0$,
			\begin{equation}\label{eq:CW}
				\Prob\big(\abs{Q(X)}\le\varepsilon\big)\le C_{\rm CW} \varepsilon^{1/2} \big(\E[Q(X)^2]\big)^{-1/4} \enspace.
			\end{equation}
			Fix $i<j$. Under $\nu$, $(u^i,u^j)\sim\mathcal N(0,I_{2d})$ and $q(z+u^i)-q(z+u^j)$ is quadratic in $(u^i,u^j)$.
			Since $q(z+u)=q(z)+2z^\top\nabla^2f(x^*)u+q(u)$, we have $q(z+u^i)-q(z+u^j)=A+B$ with $A:=q(u^i)-q(u^j)$ and $B:=2z^\top\nabla^2f(x^*)(u^i-u^j)$.
			The function $A$ is even and $B$ is odd in $(u^i,u^j)$, whose law is symmetric, so $\E[AB]=0$ and $\E[(A+B)^2]=\E[A^2]+\E[B^2]\ge\E[A^2]$.
			Moreover $q(u^i)$ and $q(u^j)$ are i.i.d.\ with
			$$
			\E[q(u^i)^2]-(\E q(u^i))^2=2\operatorname{Trace}\big(\nabla^2f(x^*)^2\big).
			$$
			Therefore
			\begin{equation*}
				\E\big[\big(q(z+u^i)-q(z+u^j)\big)^2\big]\ge\E\big[\big(q(u^i)-q(u^j)\big)^2\big]=4\operatorname{Trace}\big(\nabla^2f(x^*)^2\big) \enspace.
			\end{equation*}
			Applying \eqref{eq:CW} with $\varepsilon=\frac{16L}{3}\sigma^{1/4}$ and summing in \eqref{eq:Dsplit}, we obtain $\nu(\Dcal_\sigma(z))\le\kappa\sigma^{1/8}$ with
			\begin{equation*}
				\kappa:=\lambda2^{d/2}+\binom\lambda2C_{\rm CW}\sqrt{\frac{8L}{3}}\operatorname{Trace}\big(\nabla^2f(x^*)^2\big)^{-1/4} \enspace,
			\end{equation*}
			which is \ref{ass:coupling} with $\eta=1/8$.
		\end{proof}

		\begin{proposition}\label{prop:esA8}
			Assumption \ref{ass:domination} holds.
		\end{proposition}
		\begin{proof}
			Straightforward.
		\end{proof}

		\subsection{Application of \Cref{thm:main}}
		
		We have shown that the normalized process \eqref{eq:normalized-chain} satisfies \ref{ass:normalized}--\ref{ass:domination}. The application of \Cref{thm:main} gives the following result.

		\begin{theorem}\label{thm:es-coupled}
			Assume \ref{ass:function} and \ref{ass:quad}. Let $\beta$ be as in \eqref{eq:Vbeta} and let $\bar\sigma$ be the constant of \ref{ass:coupling} given by \Cref{prop:esA6}. 
			Then there are $\theta>0$ and $C<\infty$ such that, for every initial condition $(m_0,\sigma_0)\in\R^d\times(0,+\infty)$, there is an event $\Ecal$ with
			\begin{equation}\label{eq:B-es}
				\Prob(\Ecal^c)\le C\Big(1+\frac{\abs{m_0-x^*}^\beta}{\sigma_0^{\beta}}\Big)\sigma_0^{\theta}
			\end{equation}
			on which, for every $t\ge0$, $(m_t,\sigma_t)$ equals the state of the algorithm \eqref{eq:es} on $f_0$ from $(m_0,\sigma_0)$ with the same samples $\mathbf U$. 
			Consequently, up to a null set, on $\Ecal$
			\begin{equation}\label{eq:es-limits}
				\lim_{t\to\infty}\frac1t\log\sigma_t=\lim_{t\to\infty}\frac1t\log\big(\abs{m_t-x^*}\big)=-\CR \enspace.
			\end{equation}
			Moreover, with $z_0=(m_0-x^*)/\sigma_0$, for all $t\ge0$,
			\begin{equation}\label{eq:TV-es}
				\norm{\Lcal(z_t)-\pi_0}_\TV\le2C\Big(1+\frac{\abs{m_0-x^*}^\beta}{\sigma_0^{\beta}}\Big)\sigma_0^{\theta}+\norm{P_0^t(z_0,\cdot)-\pi_0}_\TV \enspace.
			\end{equation}
		\end{theorem}
		\begin{proof}
			By \Cref{prop:esA1,prop:esA2,prop:esA3,prop:esA4,prop:esA6,prop:esA8}, \Cref{prop:TAH} and \Cref{thm:CRquad}, Assumptions \ref{ass:normalized}--\ref{ass:domination} hold for the normalized process \eqref{eq:esmaps} with $V=V_\beta$, $\alpha=\beta/4$, $\eta=1/8$ and $\bar\sigma$ as above, and $\CR>0$. We apply \Cref{thm:main} with $z_0=(m_0-x^*)/\sigma_0$, for which $V(z_0)=1+\abs{m_0-x^*}^\beta/\sigma_0^\beta$.
			This gives the event $\Ecal$, the bound \eqref{eq:B-es} and the bound \eqref{eq:TV-es}. On $\Ecal$, by \Cref{prop:esA2}, $(m_t,\sigma_t)$ is the state of the algorithm on $f_0$ started from $(m_0,\sigma_0)$, for every $t\ge0$. The two limits in \eqref{eq:es-limits} are then those of \eqref{eq:TAH-LLN}, which hold almost surely for this algorithm, since $\CR$ is its convergence rate by \Cref{prop:esA2}.
		\end{proof}

			\subsection{Verification of \ref{assB:drift}}
				
				We show in this section that the drift assumption \ref{assB:drift} is satisfied. To that end, we describe in the next lemma the behaviour of the normalized chain \eqref{eq:normalized-chain} on linear functions.
				
			\begin{lemma}[Behaviour on linear functions]\label{lem:lin}
				Assume $\mu\le\lambda/2$ and \eqref{eq:lin-div}. Let $e\in\R^d$ with $\abs e=1$, let $f_e(x):=e^\top x$, and let $Y_{f_e}(\mathbf U_1)=\sum_{i=1}^\mu w_iU_1^{s_{f_e}(i)}$ be the mean displacement on $f_e$. Then
				\begin{equation}\label{eq:lin-facts}
					\E\big[Y_{f_e}(\mathbf U_1)\big]=-c_{\rm lin}e\qquad\text{and}\qquad\E\big[\log\Gamma(Y_{f_e}(\mathbf U_1))\big]>0 \enspace,
				\end{equation}
				where $c_{\rm lin}:=-\sum_{i\le\mu} w_i\E[\mathcal N_{i:\lambda}]>0$. Neither quantity depends on $e$.
			\end{lemma}
			\begin{proof}
				The ranking $s_{f_e}(\mathbf U_1)$ orders the coordinates $e^\top U_1^1,\dots,e^\top U_1^\lambda$ increasingly.
				Hence the vector of sorted coordinates $(e^\top U_1^{s_{f_e}(i)})_{i\le\lambda}$ has the law of $(\mathcal N_{i:\lambda})_{i\le\lambda}$.
				The components of the $U_1^i$ orthogonal to $e$ are independent of this ranking and centred.
				Hence $\E[Y_{f_e}(\mathbf U_1)]=\sum_{i\le\mu}w_i\E[\mathcal N_{i:\lambda}]e=-c_{\rm lin}e$, and $\E[\log\Gamma(Y_{f_e}(\mathbf U_1))]$ is given by \eqref{eq:CRlin}, which is positive by \eqref{eq:lin-div}.
				Neither quantity depends on $e$.
				It remains to prove $c_{\rm lin}>0$.
				By symmetry of the standard normal law, $(\mathcal N_{1:\lambda},\dots,\mathcal N_{\lambda:\lambda})$ has the law of $(-\mathcal N_{\lambda:\lambda},\dots,-\mathcal N_{1:\lambda})$, so $\E[\mathcal N_{i:\lambda}]=-\E[\mathcal N_{\lambda+1-i:\lambda}]$.
				For $i\le\mu\le\lambda/2$ we have $i<\lambda+1-i$, hence $\mathcal N_{i:\lambda}\le\mathcal N_{\lambda+1-i:\lambda}$ almost surely, with strict inequality with positive probability, so $\E[\mathcal N_{i:\lambda}]<0$.
				Since the weights are positive, $c_{\rm lin}>0$.
			\end{proof}

			We set, for constants $R_0>0$ and $r_1>0$ to be chosen below,
			\begin{equation}
				D(z) = \log\frac{\abs z}{R_0} \qquad  g(z,\sigma) = \frac{\sigma^2\abs z^2}{r_1^2}
			\end{equation}
			Then, 
			\begin{align}
				& D(z_t) - D(z_{t+1}) = \log \Gamma (Y_{t+1})  - \log \frac{\abs{z_t+c_m Y_{t+1}}}{\abs{z_t}}    \label{eq:formula-Dzt} \\
				&  g(z_t,\sigma_t)  =  \frac{\abs{m_t-x^*}^2}{r_1^2} \\
				& \sigma_t \le \frac{r_1}{R_0} e^{-D(z_t)}
			\end{align}
			
			With these definitions, we can now verify \ref{assB:drift} for the normalized chain \eqref{eq:normalized-chain}. This is done with the following lemma.
			
			\begin{lemma}\label{l:B1}
				Assume \ref{ass:function} and \ref{ass:quad}.
				There exist large $R_0>0$ and $K>0$
				and small $\bar r$, $c>0$ and $\vartheta>0$ such that for $t\in\N$, 
				if $\abs{m_t-x^*}\le\bar r$ and $\abs{z_t}\ge R_0$, then:
				\begin{enumerate}
					\item[(i)] $\E[ D(z_t) -D(z_{t+1}) | \Fcal_t ] \ge c$,
					\item[(ii)] $\E [e^{\vartheta \abs{D(z_t) -D(z_{t+1})} }| \Fcal_t ] \le K$,
					\item[(iii)] $\E[\abs{m_{t+1}-x^*}^2 | \Fcal_t  ] \le \abs{m_t-x^*}^2$.
				\end{enumerate} 
				Consequently, if $r_1\le\bar r$, then \ref{assB:drift} holds, with $\varsigma=r_1/R_0$ and $\sup_{\{D\le0\}}V_\beta=1+R_0^\beta$, since $\{D\le0\}=\{\abs z\le R_0\}$.
			\end{lemma}
			
			\begin{proof}
				Before proving the inequalities (i)-(iii), let us prove that the selection $s_{t+1}$ is close to a linear selection when $R_0$ is large.
				Indeed, on the event $(\abs{\mathbf{U}_{t+1}}_\infty\le \abs {z_t}) \supset (\abs{\mathbf{U}_{t+1}}_\infty\le R_0)$, we have
				$$
					\abs{m_t+\sigma_t {U}^i_{t+1}-x^*}
					\le  \abs{m_t-x^*}\times\left(  1+ \frac{\abs{U^i_{t+1}}}{\abs{z_t}}  \right) \le 2 \bar r \enspace.
				$$
				Consequently, by \ref{ass:function} and Taylor's inequality:
				$$
					f(m_t+\sigma_tU_{t+1}^i) = f(m_t) + \sigma_t \abs{\nabla f(m_t)} \frac{\nabla f(m_t)^\top U_{t+1}^i}{\abs{\nabla f(m_t)}} + \omega_i
				$$
				where $\abs{\omega_i} \le \frac{M}{2}\sigma_t^2\abs{\mathbf U_{t+1}}_\infty^2$ with $M:=\max_{\abs{x-x^*}\le 2\bar r} \norm{\nabla^2 f(x)}$.
				Therefore, if $s_{t+1}\neq s_{\nabla f(m_t)}$, then there exist $i<j$ such that
				$$
				\frac{\abs{\nabla f(m_t)^\top (U_{t+1}^i - U_{t+1}^j )}}{\abs{\nabla f(m_t)}}
				\le \frac{M\sigma_t\abs{\mathbf U_{t+1}}_\infty^2}{\abs{\nabla f(m_t)}}
				\le \frac{2M\abs{\mathbf U_{t+1}}_\infty^2}{ \lambda_{\min}(\nabla^2f(x^*))\abs{z_t}}
				$$
				where we have used $\abs{\nabla f(m_t)}\ge\abs{m_t-x^*} \times \lambda_{\min}(\nabla^2f(x^*))/2$.
				Thus, we have, for $e_t=\nabla f(m_t)/\abs{\nabla f(m_t)}$:
				\begin{multline*}
				\Prob[s_{t+1}\neq s_{\nabla f(m_t)} | \Fcal_t ] \\
				\le \Prob[\abs{\mathbf{U}_{t+1}}_\infty> R_0 | \Fcal_t ]
				+ \binom{\lambda}{2} \Prob \left[ \abs{e_t^\top (U_{t+1}^i - U_{t+1}^j )} \le \frac{2M\abs{\mathbf U_{t+1}}_\infty^2}{ \lambda_{\min}(\nabla^2f(x^*))R_0} \,\Big|\,\Fcal_t\right] \\
				\le \varepsilon_{R_0}
				\end{multline*}
				where $\varepsilon_{R_0}>0$ tends to $0$ when $R_0\to\infty$.
				Indeed, $e_t$ is $\Fcal_t$-measurable and $\mathbf U_{t+1}$ is independent of $\Fcal_t$, so the second conditional probability is a function of the unit vector $e_t$. By rotation invariance of the law of $\mathbf U_{t+1}$, this function is constant, and it tends to $0$ when $R_0\to\infty$ by dominated convergence, since $e^\top(U^i_{t+1}-U^j_{t+1})\ne0$ almost surely. Hence $\varepsilon_{R_0}$ depends only on $R_0$, $\bar r$, $f$, $d$ and $\lambda$.
				
				Let us prove (i).
				We have, by \eqref{eq:formula-Dzt}:
				$$
					\E[ D(z_t) -D(z_{t+1}) | \Fcal_t ] = 
					 \E[\log \Gamma (Y_{t+1})| \Fcal_t ]   - \E\left[ \log \frac{\abs{z_t+c_m Y_{t+1}}}{\abs{z_t}} | \Fcal_t \right] 
				$$
				For the first term, we have
				\begin{multline*}
					\E[\log \Gamma (Y_{t+1})| \Fcal_t ]  \ge \E[\log \Gamma (Y_{t+1}) \ind_{ s_{t+1}= s_{\nabla\! f(m_t)} } | \Fcal_t ] 
					- \E[\abs{\log \Gamma (Y_{t+1})} \ind_{ s_{t+1} \neq s_{\nabla\! f(m_t)} } | \Fcal_t ] \\
					\ge  \E[\log \Gamma (Y_{s_{\nabla\! f(m_t)}})  | \Fcal_t ]
					- 2\sqrt{\varepsilon_{R_0} \int \ell^2\dd\nu }
					\ge  \underbrace{\E[\log \Gamma (Y_{s_{\rm lin}})  ]}_{>0}	- 2\sqrt{\varepsilon_{R_0} \int \ell^2\dd\nu }
				\end{multline*}
				where we used \ref{ass:normalized} and the Cauchy--Schwarz inequality for the second to last inequality
				and \Cref{lem:lin} for the last inequality.
				For the second term, by concavity of $\log$:
				$$
					\E\left[ \log \frac{\abs{z_t+c_m Y_{t+1}}}{\abs{z_t}} | \Fcal_t \right] \le \log \left( 1 + c_m \frac{ \E[ \abs{Y_{t+1}} | \Fcal_t ]  }{\abs{z_t} }  \right) \le c_m \frac{ \E[ \max_{s\in\mathfrak{S}_\lambda} \abs{Y_{s}}  ]  }{R_0}  .
				$$
				All in all, this proves (i) given that $R_0$ is sufficiently large.
				
				We prove now (ii). We take $\vartheta\le1/8$ small enough so that $\int e^{2\vartheta\ell}\dd\nu<\infty$, which is possible by \ref{ass:domination}. Since  $\abs{\log\Gamma(Y_{t+1})}\le\ell(\mathbf U_{t+1})$ by \ref{ass:normalized}, \eqref{eq:formula-Dzt} gives
				\begin{equation*}
					e^{\vartheta \abs{D(z_t) -D(z_{t+1})} } \le e^{\vartheta\ell(\mathbf U_{t+1})} \left( \frac{\abs{z_t}^\vartheta}{\abs{z_t+c_m Y_{t+1}}^\vartheta} + \frac{\abs{z_t+c_m Y_{t+1}}^\vartheta}{\abs{z_t}^\vartheta} \right) \enspace.
				\end{equation*}
				The second term is almost surely smaller than or equal to
				\begin{equation*}
					e^{\vartheta\ell(\mathbf U_{t+1})}\Big(1+\frac{c_m}{R_0}\max_{s\in\mathfrak S_\lambda}\abs{Y_s(\mathbf U_{t+1})}\Big)^\vartheta \enspace,
				\end{equation*}
				whose expectation is finite and does not depend on $(z_t,\sigma_t)$. For the first term, by Cauchy--Schwarz:
				\begin{equation*}
					\E \left[ e^{\vartheta\ell(\mathbf U_{t+1})} \frac{\abs{z_t}^\vartheta}{\abs{z_t+c_m Y_{t+1}}^\vartheta} \Big| \Fcal_t \right] 
					\le \left( \int e^{2\vartheta \ell} \dd\nu \right)^{1/2} \left( \E \left[ \frac{\abs{z_t}^{2\vartheta}}{\abs{z_t+c_m Y_{t+1}}^{2\vartheta}} \Big| \Fcal_t \right] \right)^{1/2} \enspace.
				\end{equation*}
				On the event $\{c_m\abs{Y_{t+1}}\le\abs{z_t}/2\}\subset\{\abs{z_t+c_mY_{t+1}}\ge\abs{z_t}/2\}$,  ${\abs{z_t}^{2\vartheta}}/{\abs{z_t+c_m Y_{t+1}}^{2\vartheta}} \le 2^{2\vartheta}$.
				Besides, by Cauchy--Schwarz:
				\begin{multline*}
					\E \left[ \frac{\abs{z_t}^{2\vartheta}}{\abs{z_t+c_m Y_{t+1}}^{2\vartheta}} \ind_{\{c_m\abs{Y_{t+1}}>\abs{z_t}/2\}} \Big| \Fcal_t \right] \\
					\le 
					\underbrace{\abs{z_t}^{2\vartheta} \Prob\Big( c_m\max_{s\in\mathfrak S_\lambda}\abs{Y_s(\mathbf U_{t+1})} > \frac{\abs{z_t}}{2} \Big)^{1/2}}_{\text{bounded when } \abs{z_t}>R_0}
					 \left( \E \left[ \abs{z_t+c_m Y_{t+1}}^{-4\vartheta} \Big| \Fcal_t \right] \right)^{1/2} \enspace.
				\end{multline*}
				 Since $4\vartheta\le1/2$, we have
				\begin{align*}
					\E \left[ \abs{z_t+c_m Y_{t+1}}^{-4\vartheta} \Big| \Fcal_t \right]
					&\le \sum_{s\in\mathfrak S_\lambda} \E \left[ \max\big(1,\abs{z_t+c_m Y_s(\mathbf U_{t+1})}^{-1/2}\big) \Big| \Fcal_t \right] \\
					&\le \lambda! \left( 1 + (2\pi c_m^2\abs w^2)^{-d/2} \int_{\abs y\le1} \abs y^{-1/2} \dd y \right) \enspace,
				\end{align*}
				which gives (ii). For (iii), we note that
				\begin{align*}
					\E[\abs{m_{t+1}-x^*}^2 | \Fcal_t  ] 
					& = \abs{m_t-x^*}^2 + 2 c_m \sigma_t \E[ Y_{t+1}^\top (m_t-x^*) | \Fcal_t  ] + c_m^2\sigma_t^2 \E[ \abs{Y_{t+1}}^2  | \Fcal_t  ] \\
					& \le \abs{m_t-x^*}^2 + \frac{\abs{m_t-x^*}^2}{\abs{z_t}} \times \Big(   2c_m \E\Big[ \frac{Y_{t+1}^\top (m_t-x^*)}{\abs{m_t-x^*}} \,\Big|\, \Fcal_t  \Big] \\
					&\qquad + c_m^2 \frac{\E[ \max_{s\in\mathfrak{S}_\lambda} \abs{Y_{s}}^2 ] }{\abs{z_t}}  \Big).
				\end{align*}
				However,
				\begin{multline*}
					\E\left[ \frac{Y_{t+1}^\top (m_t-x^*)}{\abs{m_t-x^*}} \,\Big|\, \Fcal_t  \right]
					= \E\left[ \frac{Y_{t+1}^\top (m_t-x^*)}{\abs{m_t-x^*}}  \ind_{\{ s_{t+1}= s_{\nabla f(m_t)}\}} \,\Big|\, \Fcal_t  \right] \\
					+  \E\left[ \frac{Y_{t+1}^\top (m_t-x^*)}{\abs{m_t-x^*}}  \ind_{\{ s_{t+1}\neq s_{\nabla f(m_t)}\}} \,\Big|\, \Fcal_t  \right] \\
					\le \E\left[ \frac{Y_{s_{\nabla f(m_t)}}^\top (m_t-x^*)}{\abs{m_t-x^*}}  | \Fcal_t  \right]
					+ 2\sqrt{\varepsilon_{R_0} \E\left[\max_{s\in\mathfrak{S}_\lambda} \abs{Y_s}^2 \right]   }
				\end{multline*}
				by the Cauchy--Schwarz inequality as in the proof of (i), for which
				$$
				\E\left[ \frac{Y_{s_{\nabla f(m_t)}}^\top (m_t-x^*)}{\abs{m_t-x^*}}  | \Fcal_t  \right] 
				= -c_{\rm lin} \frac{\nabla f(m_t)^\top (m_t-x^*)}{\abs{\nabla f(m_t)}\abs{m_t-x^*}}
				<0
				$$ 
				by \Cref{lem:lin}.
				Besides, when $\abs{m_t-x^*}\le \bar r$, the Taylor identity
				$$
					\nabla f(m_t) = \underbrace{\nabla f(x^*)}_{=0} + \nabla^2 f(x^*) (m_t-x^*) + h_{m_t-x^*} \times \abs{m_t-x^*}
				$$
				with $h_{m_t-x^*} \to 0$ when $m_t\to x^*$, gives
				for sufficiently small $\bar r>0$:
				\begin{gather*}
				\nabla f(m_t)^\top (m_t-x^*) \ge \frac{\lambda_{\min}(\nabla^2 f(x^*))}{2} \abs{m_t-x^*}^2, \\
				\abs{\nabla f(m_t)} \le 2\lambda_{\max}(\nabla^2 f(x^*)) \abs{m_t-x^*},
				\end{gather*}
				thus
				$$
				\E\left[ \frac{Y_{s_{\nabla f(m_t)}}^\top (m_t-x^*)}{\abs{m_t-x^*}}  | \Fcal_t  \right] 
				\le -c_{\rm lin}
				\frac{ \lambda_{\min}(\nabla^2 f(x^*)) }{4\lambda_{\max}(\nabla^2 f(x^*)) } < 0 .
				$$
				Therefore, taking $R_0$ sufficiently large yields
				$$
				 2c_m \E\left[ \frac{Y_{t+1}^\top (m_t-x^*)}{\abs{m_t-x^*}} | \Fcal_t  \right] + c_m^2 \frac{\E[ \max_{s\in\mathfrak{S}_\lambda} \abs{Y_{s}}^2 ] }{\abs{z_t}} \le 0 ,
				$$
				which finally proves (iii).
			\end{proof}
			
			\subsection{Proof of \Cref{thm:es2}}
			
			We are able now to conclude the proof of \Cref{thm:es2}.
			
			By \Cref{prop:esA1,prop:esA2,prop:esA3,prop:esA4,prop:esA6,prop:esA8}, \Cref{prop:TAH}, \Cref{thm:CRquad} and \Cref{l:B1}, \ref{ass:normalized}--\ref{ass:domination} and \ref{assB:drift} hold for every $r_1\in(0,\bar r]$.
			Therefore, by \Cref{th:decoupled-initial-conditions}, the event $\Ecal_\tau$ exists and
			\begin{equation}\label{eq:es-bound}
				\Prob(\Ecal_\tau^c) \le \frac{\abs{m_0-x^*}^2}{r_1^2} + C (1+R_0^\beta) \times \left( \sigma_0^\theta+ C_{\rm occ} \times \left( \frac{r_1}{R_0} \right)^\theta \int e^{\theta \ell}\dd\nu  \right) \enspace.
			\end{equation}
			The constants $C$ and $\theta$ come from \Cref{thm:main} and depend only on the constants of \ref{ass:normalized}--\ref{ass:domination}, and $C_{\rm occ}$ depends only on $(c,\vartheta,\theta,K)$.
			Let $\varepsilon\in(0,1)$. We choose first $r_1\in(0,\bar r]$ such that the last term of \eqref{eq:es-bound} is at most $\varepsilon/3$, then $s_\varepsilon>0$ such that $C(1+R_0^\beta)s_\varepsilon^\theta\le\varepsilon/3$, and finally $r_\varepsilon:=r_1\sqrt{\varepsilon/3}$.
			For $\abs{m_0-x^*}\le r_\varepsilon$ and $\sigma_0\le s_\varepsilon$, \eqref{eq:es-bound} gives $\Prob(\Ecal_\tau^c)\le\varepsilon$.

			On $\Ecal_\tau$, for $t\ge\tau$, $(z_t,\sigma_t)$ matches the limit chain \eqref{eq:limit-chain} started at time $\tau$ from $(z_\tau,\sigma_\tau)$. By \Cref{prop:esA2}, this is the normalized process and the step size of the algorithm \eqref{eq:es} on $f_0$ started from $(m_\tau,\sigma_\tau)$, with $m_t-x^*=\sigma_tz_t$ for every $t$.
			Since $\tau$ is an almost surely finite stopping time, the samples $\mathbf U_{\tau+1},\mathbf U_{\tau+2},\dots$ are i.i.d.\ with law $\nu$ and independent of $\Fcal_\tau$, so \eqref{eq:TAH-LLN} applied from time $\tau$ gives the two limits of \Cref{thm:es2} almost surely on $\Ecal_\tau$.

	\section{Proof of \Cref{thm:CRquad}}\label{sec:CRquad}
	
	The proof of the local convergence result in \Cref{thm:es2} relies on the positivity of the convergence rate on quadratic functions, as stated in \Cref{thm:CRquad}. The present section is devoted to prove that property.

	Let $f(x)=(x-x^*)^\top H(x-x^*)$ with $H$ symmetric positive definite, and assume \eqref{eq:lin-div}, nonincreasing weights and $\mu\le\lambda/2$, as in \Cref{thm:CRquad}.
	With the notation of \Cref{sec:es-notation}, let $Y_z(u):=Y_{f(x^*+z+\cdot)}(u)$.
	By \eqref{eq:zrec} the normalized process is the Markov chain
	\begin{equation}\label{eq:chain}
		z_{t+1}=\frac{z_t+c_mY_{z_t}(\mathbf U_{t+1})}{\Gamma(Y_{z_t}(\mathbf U_{t+1}))}
	\end{equation}
	with kernel $P$. By \Cref{prop:TAH}, $P$ is positive recurrent with a unique invariant probability $\pi$.
	We write $Y_z=Y_z(\mathbf U_1)$. By \eqref{eq:CSA2} and \eqref{eq:CRdef},
	\begin{equation}\label{eq:g}
		\bar\gamma(z):=
		2 d_\sigma \E\log\Gamma(Y_z)
		=\frac{\E\abs{Y_z}^2}{d\abs w^2}-1\ge-1,\qquad \CR=-\frac{\pi(\bar\gamma)}{2d_\sigma} \enspace.
	\end{equation}
	By \Cref{lem:lin}, for a unit vector $e$, $\E[Y_{f_e}(\mathbf U_1)]=-c_{\rm lin}e$ with $c_{\rm lin}>0$.
	Also by \Cref{lem:lin}, $\E\log\Gamma(Y_{f_e}(\mathbf U_1))>0$, that is, $\E\abs{Y_{f_e}(\mathbf U_1)}^2>d\abs w^2$.

	We use the bounds $\abs{Y_z(u)}\le C\abs u$ and $\abs{\log\Gamma(Y_z(u))}\le C(1+\abs u^2)/d_\sigma$ for some well-chosen constant $C>0$.
	We take $d_\sigma$ large enough so that
	\begin{equation}\label{eq:ell}
		2\abs{\log\Gamma(Y_z(u))}\le\frac{1+\abs u^2}{4}\qquad\text{for all }z\text{ and }u \enspace.
	\end{equation}

	Besides, the bound $\abs{Y_z(u)}\le C\abs u$ gives $\abs{\bar\gamma}\le C$.
	Note that the probability that $Y_{z'}(u)=Y_z(u)$ tends to one when $z'$ tends to $z$. By the dominated convergence theorem, $\bar\gamma$ and $z\mapsto\E[z^\top HY_z]$ are continuous.
	By \cite[Lemma 5.6]{gissler2024linear}, item (i), applied with an identity covariance matrix, and with $\mu_{\rm eff}=1/\abs w^2$ there, $\E\abs{Y_0}^2<d\abs w^2$ for nonincreasing weights, i.e.,
	\begin{equation}\label{eq:g0}
		\bar\gamma(0)<0 \enspace.
	\end{equation}
	By rotation invariance of the algorithm and of the law of $\mathbf U_1$, we may assume that $H$ is diagonal.
	Then \cite[Proposition 4.17]{gissler2024linear}, applied with an identity covariance matrix, gives $\E[Y_z]_j<0$ for every $j$ with $z_j>0$ and nonincreasing weights, hence $\E[Y_z]_j>0$ when $z_j<0$ and $\E[Y_z]_j=0$ when $z_j=0$ by the symmetry $u_j\mapsto-u_j$.
	Since $H$ is diagonal with positive entries, $z^\top H\E[Y_z]=\sum_jH_{jj}z_j\E[Y_z]_j$ and for $z\ne0$,
	\begin{equation}\label{eq:drift-sign}
		\E[z^\top HY_z]<0\qquad\text{for every }z\ne0 \enspace.
	\end{equation}
	As in the first step of the proof of \Cref{l:B1}, for $z\ne0$ the rankings by $f(x^*+z+\cdot)$ and by the linear function $f_{Hz/\abs{Hz}}$ coincide outside an event whose probability depends only on $\abs z$ and tends to $0$ as $\abs z\to\infty$.
	Since both displacements are bounded by $C\abs{\mathbf U_1}$, the Cauchy--Schwarz inequality gives
	\begin{equation}\label{eq:far}
		\lim_{\abs z\to\infty}\bar\gamma(z)=\frac{\E\abs{Y_{f_e}(\mathbf U_1)}^2}{d\abs w^2}-1>0\qquad\text{and}\qquad\lim_{\abs z\to\infty}\frac{\E[z^\top HY_z]}{\abs{Hz}}=-c_{\rm lin}<0 \enspace.
	\end{equation}
	Finally, let
	\begin{equation}\label{eq:W}
		W(z):=\sqrt{1+z^\top Hz} \enspace.
	\end{equation}
	
	The next lemma shows a Foster--Lyapunov condition required to apply a comparison theorem~\cite[Theorem 14.3.7]{meyn2009markov}.
	\begin{lemma}\label{lem:drift}
		If $c_m\le1$ and $d_\sigma$ is large enough, for all $z\in\R^d$,
		\begin{equation}\label{eq:drift}
			PW(z)-W(z)\le c_m\frac{\E[z^\top HY_z]}{W(z)}+\frac{C}{d_\sigma}+Cc_m^2 \enspace.
		\end{equation}
	\end{lemma}
	\begin{proof}
		Fix $z$. Since $z_1=(z+c_mY_z)/\Gamma(Y_z)$ when $z_0=z$, see \eqref{eq:chain}, then,
		\begin{equation*}
			W(z_1)^2-W(z)^2=\Gamma(Y_z)^{-2}\big(z^\top Hz+2c_mz^\top HY_z+c_m^2Y_z^\top HY_z\big)-z^\top Hz \enspace.
		\end{equation*}
		By concavity of the square root,
		\begin{multline}\label{eq:split}
			W(z_1)-W(z)\le\frac{W(z_1)^2-W(z)^2}{2W(z)}\\
			=\frac{(\Gamma(Y_z)^{-2}-1)z^\top Hz}{2W(z)}+\frac{c_m\Gamma(Y_z)^{-2}z^\top HY_z}{W(z)}+\frac{c_m^2\Gamma(Y_z)^{-2}Y_z^\top HY_z}{2W(z)} \enspace.
		\end{multline}
		We take expectations and bound the three terms, using \eqref{eq:ell} and $\E e^{\abs{\mathbf U_1}^2/4}<\infty$.

		For the last term, for a sufficiently large constant $C>0$ (independent of $z$, $c_m$ and $d_\sigma$), we have
		 $$
		 \frac{c_m^2 \E[\Gamma(Y_z)^{-2}Y_z^\top HY_z]}{2W(z)} \le
		 Cc_m^2 \enspace,
		 $$ 
		 since $\Gamma(Y_z)^{-2}\le e^{1/d_\sigma}\le e$ and $W\ge1$.

		For the middle term, $\abs{z^\top HY_z}\le CW(z)\abs{\mathbf U_1}$ and
		\begin{equation*}
			\abs{\Gamma(Y_z)^{-2}-1}\le2\abs{\log\Gamma(Y_z)}e^{2\abs{\log\Gamma(Y_z)}}\le \frac{C}{d_\sigma}(1+\abs{\mathbf U_1}^2)e^{(1+\abs{\mathbf U_1}^2)/4} \enspace,
		\end{equation*}
		so
		\begin{equation*}
			\E\Big[\frac{c_m\Gamma(Y_z)^{-2}z^\top HY_z}{W(z)}\Big]\le c_m\frac{\E[z^\top HY_z]}{W(z)}+\frac{Cc_m}{d_\sigma} \enspace.
		\end{equation*}

		For the first term, by Taylor's inequality,
		\begin{equation*}
			\Gamma(Y_z)^{-2}\le1-2\log\Gamma(Y_z)+2(\log\Gamma(Y_z))^2e^{2\abs{\log\Gamma(Y_z)}} \enspace.
		\end{equation*}
		Moreover, \eqref{eq:g} gives $\E\Gamma(Y_z)^{-2}-1\le-\bar\gamma(z)/d_\sigma+C/d_\sigma^2$.
		By \eqref{eq:far}, up to taking $C$ larger, we have $\bar\gamma\ge1/C$ on $\{\abs z>C\}$,
		so, for $d_\sigma\ge C^2$,
		$$
		\frac{(\E\Gamma(Y_z)^{-2}-1)z^\top Hz}{2W(z)} \le 0 \qquad \text{for } \abs z >C \enspace.
		$$
		 On $\{\abs z\le C\}$ we have $z^\top Hz/(2W(z))\le C$ and $-\bar\gamma(z)\le1$, hence
		 $$
		 \frac{(\E\Gamma(Y_z)^{-2}-1)z^\top Hz}{2W(z)} \le
		 C/d_\sigma \enspace.
		 $$
		Summing the three bounds and using $c_m\le1$ in the middle one gives \eqref{eq:drift}.
	\end{proof}
	
	We conclude now the proof of \Cref{thm:CRquad}.

	\begin{proof}[Proof of \Cref{thm:CRquad}]
		By \eqref{eq:g0} and the continuity of $\bar\gamma$, the set $\{\bar\gamma>\bar\gamma(0)/2\}$ is at a positive distance from $0$.
		By \eqref{eq:drift-sign}, the continuous function $z\mapsto-\E[z^\top HY_z]/W(z)$ is positive on $\R^d\setminus\{0\}$. For $\abs z\ge1$ it is the product of $-\E[z^\top HY_z]/\abs{Hz}$ and $\abs{Hz}/W(z)\ge1/C$, so its $\liminf$ as $\abs z\to\infty$ is positive by \eqref{eq:far}. Hence, for $C>0$ sufficiently large,
		\begin{equation}\label{eq:inf}
			-\frac{\E[z^\top HY_z]}{W(z)}\ge\frac1C\qquad\text{on }\{\bar\gamma>\bar\gamma(0)/2\} \enspace.
		\end{equation}

		The functions $W$ and $z\mapsto-c_m\E[z^\top HY_z]/W(z)$ and the constant in \eqref{eq:drift} are nonnegative.
		By \Cref{lem:drift},
		 \cite[Theorem 14.3.7]{meyn2009markov} applies and gives
		\begin{equation*}
			\pi(\bar\gamma>\bar\gamma(0)/2)\le C\int\frac{-\E[z^\top HY_z]}{W(z)}\pi(\dd z)\le C\Big(c_m+\frac{1}{c_md_\sigma}\Big) \enspace.
		\end{equation*}
		Since $\abs{\bar\gamma}\le C$ everywhere, we have $\bar\gamma-\bar\gamma(0)/2\le2C$ on $\{\bar\gamma>\bar\gamma(0)/2\}$ and $\bar\gamma-\bar\gamma(0)/2\le0$ elsewhere, so
		\begin{equation*}
			\pi(\bar\gamma)\le\frac{\bar\gamma(0)}{2}+2C\pi(\bar\gamma>\bar\gamma(0)/2)\le\frac{\bar\gamma(0)}{2}+2C\Big(c_m+\frac{1}{c_md_\sigma}\Big)\le\frac{\bar\gamma(0)}{4}
		\end{equation*}
		if $c_m\le-\bar\gamma(0)/(16C)$ and $c_md_\sigma\ge-16C/\bar\gamma(0)$.
		These two conditions give $c_m\le1$ and $d_\sigma\ge(16C/\bar\gamma(0))^2$, which make $c_m$ and $d_\sigma$ as required in \eqref{eq:ell} and in \Cref{lem:drift}, up to taking $C$ larger.
		 By \eqref{eq:g}, $\CR\ge-\bar\gamma(0)/(8d_\sigma)>0$, which ends the proof.
	\end{proof}

	\bibliographystyle{alpha}
	\bibliography{references}
	
\end{document}